\documentclass[a4paper,11pt]{article}
\usepackage[utf8]{inputenc}
\usepackage[english]{babel}

\usepackage{amsmath,amssymb,amsthm}
\usepackage{booktabs}
\usepackage{fullpage}
\usepackage[hidelinks]{hyperref}

\usepackage{cleveref}
\usepackage{dsfont}
\usepackage{enumerate}
\usepackage[shortlabels]{enumitem}
\usepackage{float}
\usepackage{makecell}
\usepackage{orcidlink}
\usepackage{tabularx}
\usepackage{multirow}
\usepackage{multicol}
\usepackage{natbib}
\usepackage{tikz}
\usetikzlibrary{shapes.geometric}
\usetikzlibrary{positioning}
\usetikzlibrary{arrows,decorations.markings}
\providecommand{\keywords}[1]
{
    {\noindent\footnotesize\textbf{Keywords---} #1}
}
\providecommand{\acknowledgement}[1]
{
    {\noindent\textbf{Acknowledgements---} #1}
}

\newtheorem{theorem}{Theorem}[section]

\newtheorem{lemma}[theorem]{Lemma}
\newtheorem{proposition}[theorem]{Proposition}
\newtheorem{corollary}[theorem]{Corollary}
\theoremstyle{definition}

\newtheorem{definition}[theorem]{Definition}

\theoremstyle{remark}

\newcommand{\abs}[1]{\left\lvert#1\right\rvert}
\newcommand{\norm}[1]{\left\lVert#1\right\rVert}
\newcommand{\enc}[1]{\left(#1\right)}
\newcommand{\setEnc}[1]{\left\{#1\right\}}
\newcommand{\defeq}{\mathrel{\rlap{\raisebox{0.3ex}{$\cdot$}}\raisebox{-0.3ex}{$\cdot$}}=}
\newcommand{\smid}{\,\middle\vert\,}
\newcommand{\blankArg}{{{}\cdot{}}}
\newcommand{\approxOp}[1]{\hat{#1}}

\newcommand{\nat}{\mathds{N}}

\newcommand{\arc}{a}
\newcommand{\arcs}{A}
\newcommand{\block}{B}
\newcommand{\cost}{c}
\newcommand{\course}{C}
\newcommand{\dep}{d}
\newcommand{\deps}{\mathcal{D}}
\newcommand{\consumption}{e}
\newcommand{\minReqChrg}{E}
\newcommand{\ODEfunc}{f}
\newcommand{\curves}{F}
\newcommand{\DAG}{G}

\newcommand{\commodity}{k}
\newcommand{\commodities}{K}
\newcommand{\mixLower}{\ell}
\newcommand{\numLinSegments}{m}
\newcommand{\node}{n}
\newcommand{\nodes}{N}
\newcommand{\arcPath}{p}
\newcommand{\arcPaths}{P}

\newcommand{\chrgslot}{s}
\newcommand{\chrgslots}{\mathcal{S}}
\newcommand{\timevar}{t}
\newcommand{\mixUpper}{u}
\newcommand{\vcl}{v}
\newcommand{\vcls}{\mathcal{V}}
\newcommand{\voltage}{V}
\newcommand{\flowvar}{x}
\newcommand{\chrgvar}{y}
\newcommand{\maxChrg}{Y}
\newcommand{\slope}{\alpha}
\newcommand{\offset}{\beta}
\newcommand{\outset}{\delta^+}
\newcommand{\inset}{\delta^-}
\newcommand{\diffOp}[1]{\Delta#1}
\newcommand{\errorFunc}{\varepsilon}
\newcommand{\chrgcurveMax}{\zeta}
\newcommand{\horizonEnd}{\mathcal{H}}
\newcommand{\timestep}{\theta}
\newcommand{\timeOp}[1]{\Theta#1}
\newcommand{\mixCoeff}{\kappa}
\newcommand{\mixTypes}{\mathcal{K}}
\newcommand{\mixCons}{\bar{K}}
\newcommand{\chrgcurve}{\xi}
\newcommand{\numChrgEvents}{\sigma}
\newcommand{\trip}{\tau}
\newcommand{\trips}{\mathcal{T}}
\newcommand{\chrgdiffvar}{\varphi}
\newcommand{\chrgdiffDomain}{\Phi}
\newcommand{\linSpline}[1]{\Psi#1}
\newcommand{\maxPowerDraw}{\Omega}
\newcommand{\maxPowerCoeff}{\omega}

\newcommand{\accessPoint}{\tilde{\chrgslot}}
\newcommand{\gridAccess}{\tilde{\chrgslots}}

\newcommand{\ODEfuncCC}{\ODEfunc_{CC}}
\newcommand{\ODEfuncCV}{\ODEfunc_{CV}}

\newcommand{\CCCVbreakTime}{\timevar_{\voltage}}
\newcommand{\chrgtimeMax}{\timevar_{full}}
\newcommand{\chrgStart}{\timevar_s}
\newcommand{\chrgEnd}{\timevar_e}

\newcommand{\CCCVbreakSoC}{\chrgvar_{\voltage}}
\newcommand{\initSoC}{\chrgvar_s}
\newcommand{\finalSoC}{\chrgvar_e}

\newcommand{\commoditiesE}{\commodities^{E}}
\newcommand{\chrgNodes}{\bar{\chrgslots}}
\newcommand{\depotPullouts}{\arcs_{\deps}^+}
\newcommand{\rechargeArcs}{\arcs^E}
\newcommand{\costtou}{\bar{\cost}}
\newcommand{\inPaths}{\arcPaths^-}
\newcommand{\outPaths}{\arcPaths^+}

\newcommand{\VSP}{\emph{VSP}}
\newcommand{\BSP}{\emph{BSP}}
\newcommand{\EVSP}{\emph{EVSP}}
\newcommand{\EBSP}{\emph{EBSP}}
\newcommand{\FRVCP}{\emph{FRVCP}}
\newcommand{\soc}{\emph{soc}}
\newcommand{\tou}{\emph{tou}}
\newcommand{\cccv}{\emph{CC-CV}}

\newcommand{\features}{\ref{feat:nonlinear}-\ref{feat:mixed}}
\newcommand{\rmlinear}{L}
\newcommand{\rmexpansion}{ESE}
\newcommand{\rmspline}{LS}
\newcommand{\rmalgo}{EiA}
\newcommand{\rmrate}{CID}

\newcolumntype{Y}{>{\raggedleft\arraybackslash}X}

\newcommand{\myendproof}{\hfill\mbox{$\square$}}
\newcommand{\myendclaimproof}{\hfill\mbox{$\diamond$}}

\title{Electric Bus Scheduling with Non-Linear Charging, Power Grid Bottlenecks, and Dynamic Recharge Rates}
\author{%
\begin{tabular}{ccc}
Fabian L{\"o}bel \orcidlink{0000-0001-5433-184X} & Ralf Bornd{\"o}rfer \orcidlink{0000-0001-7223-9174} & Steffen Weider \\
Zuse Institute Berlin & Zuse Institute Berlin & IVU Traffic Technologies AG \\
\texttt{fabian.loebel@zib.de} & \texttt{borndoerfer@zib.de} & \texttt{stw@ivu.de}
\end{tabular}%
}
\date{\today}

\begin{document}

\maketitle

\begin{abstract}
Public transport operators are gradually electrifying their bus fleets, predominantly with battery-powered drive trains.
These buses commonly have to be scheduled to recharge in-service, which gives rise to a number of challenges.
A major problem is that the relationship between charging time and replenished driving range is non-linear, which is often approximately modeled.
We examine the associated approximation error and show how it can result in a gross over- or underestimation of the fleet size.
Moreover, we demonstrate that commonly used piecewise linear underestimations of the charge curve do not result in an underestimation of the predicted charge states in electric vehicle scheduling and routing models.

Furthermore, since power grid upgrades are currently not keeping up with an ever growing electricity demand, operators are introducing active charge management tools to dynamically adjust the charging speed depending on the amount of available energy.
It is therefore imperative to extend electric bus scheduling models to account for these developments.

We propose a novel mixed-integer programming formulation for the electric bus scheduling problem featuring an improved approximation of the non-linear battery charging behavior as well as dynamic recharge speeds to accommodate grid load limits.
The idea is to linearly interpolate what we call the charge increment function, which is closely related to the derivative of the commonly used charge curve.
This provides very good error control and integrates easily into integer programming models. 
We demonstrate the practical usefulness of our model on a diverse library of real-life instances.
\end{abstract}


\section{Introduction} \label{sec:intro}
Public transport operators are gradually shifting their procurement strategy for their bus fleets from combustion to electric drive trains.
In Germany, for example, only $55\%$ of new acquisitions may be powered by fossil fuels, with the quota decreasing to $35\%$ in 2026 \citep{SL23}.
Major cities like Berlin, Hamburg, or Munich have pledged to fully electrify their public transport systems by 2030 \citep{Neissen23}.
The electric bus market is dominated by battery-powered vehicles, which operators prefer to charge at their depots and selected terminals with fast chargers to minimize acquisition costs as well as political and legal concerns \citep{SGCL19,EC19,JG20,Nuremberg23}.

It is well known that electric battery-powered buses (\emph{ebuses}) have shorter driving ranges and longer refueling duration than their combustion counterparts.
This often requires recharge events in-service to maximize the available range throughout the day, especially during winter and summer when environmental conditions hamper available ranges even further \citep[cf.][]{SGCL19}.
For example, air conditioning may reduce driving ranges of 250 km to 175 km in practice \citep{WLLJ22}.
Therefore, electric vehicle schedules, i.e., sequences of passenger trips serviced by the same bus, have to include downtime and outright detours to charging infrastructure to be able to operate the same networks that combustion powered buses can on a single refuel action upon returning to the depot at the end of the day.
The relationship between charge duration and state of charge (\soc) is in general non-linear \citep[cf.][]{PJLV17}, which should be reflected in electric vehicle scheduling models (c.f. \cite{MGMV17}, \cite{OK20}, \cite{ZLTDCGL21}, as well as \Cref{sec:discParam}).

Moreover, operators like Hamburg's \emph{Hochbahn AG} and \emph{VHH} \citep{JES19} have adopted active charge management systems to take advantage of fluctuating energy prices throughout the day.
These systems monitor all currently plugged in buses and charge them such that the energy costs are optimized while minimum charge levels are maintained for flexible vehicle disposition \citep{Hamburg23}.
They accomplish this by throttling or increasing the charging rate as needed.
There are even case studies where rural operators rent out their idle bus batteries as grid buffers outside of the peak school traffic hours \citep{Leptien24}, achieving net negative electricity costs by discharging into the power grid \citep{Losel23}.

Furthermore, electricity providers are imposing load-dependent energy prices and hard grid load limits restricting the total amount of power an operator can extract from the grid to fuel their fleet at any given moment \citep[cf.][]{WLLJ22}, because the capacity of the power grid is becoming a major bottleneck for the energy transition \citep{IEA23}. 
This is felt especially by companies whose electricity demand is going from residential to industrial levels like public transport operators, whose growing electric fleets strain the power grid further than the currently installed rating \citep[e.g.][]{JES19,SGCL19,ABBCK20,BPLP20,WLLJ22,Neissen23}.
Some public transport operators in Germany have already taken steps to upgrade the connected wattage at their vehicle facilities if possible like in Hamburg \citep{RS23} and D{\"u}sseldorf \citep{Richarz24}, or produce electricity on-site from wind or solar sources like in Nuremberg \citep{Nuremberg23}, but these transformations are costly as they often involve having to move facilities to entirely new locations.
This further increases the need for active charge management systems to achieve so-called peak shaving, where the grid load is distributed as evenly about the day as possible \citep{AR24}.
Of course, this is subject to the time windows allowed by the underlying vehicle schedule, and should therefore be included in electric bus scheduling models by keeping recharge speeds dynamic.

\section{The Electric Bus Scheduling Problem} \label{sec:problemdef}
The literature knows two variants of the vehicle scheduling problem.
In general, one arises from the public transport sector and the other from more general logistic and distribution settings.
\begin{definition}
    Given a set of vehicle types $\vcls$, depots $\deps$, timetabled vehicle duties $\trips$ called \emph{trips} and knowledge about feasible or admissible connections, the \emph{Bus Scheduling Problem} (\BSP) is to partition the trips into \emph{courses} such that
    \begin{itemize}
        \item every course is assigned to a bus type and depot,
        \item the bus of each course can service its trips in order and return to the depot of origin,
        \item and the costs are minimized.
    \end{itemize}
\end{definition}
$\trips$ would usually consist of the passenger trips arising from having to service periodically timetabled lines.
\BSP{} can be modeled using an acyclic graph on which a minimum cost binary flow hitting every trip node has to be found.
Alternatively, one can view $\trips$ and the admissible connections as a partially ordered set and the \BSP{} as a chain partition problem.
The \BSP{} is polynomially solvable if $\abs{\vcls} = \abs{\deps} = 1$ and no additional side constraints are imposed, otherwise it is NP-hard in general \citep{BCG87,BK10}.

The \emph{Vehicle Scheduling Problem} (\VSP) generalizes the \BSP{} by considering customers instead of vehicle duties or trips.
Vehicles have to deliver commodities to the customers, usually within certain time windows, satisfying a priori known demands and load capacities per vehicle type.
Both problem variants have been studied for decades.
However, fleet electrification poses new challenges for modeling and solving approaches.
\begin{definition}
    Given battery capacities per vehicle type, energy consumptions on relevant activities, and a set of charging slots $\chrgslots$, the \emph{Electric Bus Scheduling Problem} (\EBSP) extends the \BSP{} by requiring that no battery is fully depleted throughout the service horizon while permitting detours to chargers to replenish driving range.
\end{definition}
We call a sequence of trips that can be serviced in order without interruptions a \emph{block}.
Courses in the electric setting therefore become alternating sequences of blocks and recharge events.
The \emph{Electric Vehicle Scheduling Problem} (\EVSP) extends the \VSP{} analogously.
We explicitly focus on \EBSP{} in this paper.

From the operator requirements outlined earlier and gaps in the literature, see \Cref{sec:literaturReview}, we conclude that an \EBSP{} model should include the following features:
\begin{enumerate}[(i)]
    \item \label{feat:nonlinear} non-linear charging behavior,
    \item \label{feat:partial} partial charging (i.e., ebuses may depart charge stations before they reach full charge),
    \item \label{feat:rates} dynamic recharge rates (i.e., ebuses do not have to be charged at full power),
    \item \label{feat:ToU} time-of-use (\tou) electricity prices,
    \item \label{feat:capacity} charger slot capacities,
    \item \label{feat:gridload} grid load limits,
    \item \label{feat:mixed} and mixed fleets of electric and non-electric types.
\end{enumerate}
Feature \ref{feat:mixed} is important for the current transition period as fleets are not going to be fully electrified instantaneously.

The remainder of this paper is structured as follows:
In \Cref{sec:literaturReview} we review related literature.
In \Cref{sec:rechargeModeling} we examine the numerical impact of approximating the charging behavior, in particular concerning worst-case instances and state-of-the-art practices.
In \Cref{sec:CID} we extend upon a charge curve linearization scheme we introduced in \citep{LBWatmos23}.
In \Cref{sec:MIP} we accordingly propose a novel \EBSP{} mixed-integer program fulfilling \features{}.
In \Cref{sec:study} we conclude with a computational study on a set of diverse but anonymous real-life \EBSP{} instances.

\section{Related Literature} \label{sec:literaturReview}
\cite{EC19} extensively review 75 contributions on the \EVSP{} and the related electric vehicle routing problem up until early 2019.
They find that the literature has been neglecting non-linear battery behavior, robustness of solutions, and that there is a tendency to assume an infinite number of charger slots at recharge facilities.
Fixing charge rates a priori has not been addressed by \citeauthor{EC19}, but from our own review, none of the papers featured concern themselves with it, either.
The most popular solving approaches are genetic algorithms and neighborhood searches.

\cite{PLL22} review 23 contributions on \EBSP{} until 2021.
They have four papers in common with the survey by \citeauthor{EC19}.
\citeauthor{PLL22} arrive at similar conclusions as \citeauthor{EC19} and additionally point out that few papers even consider partial charging, i.e., most papers require ebuses to always recharge fully upon visiting a charging facility.
For solving the \EBSP{}, linear programming approaches, in particular column generation, are more common, comprising roughly half the presented literature.

\citeauthor{PLL22} further survey publications on \emph{charging scheduling} where a fixed bus schedule is part of the input and the goal is to find the most efficient recharge policy compliant with the given vehicle downtime.
Since the problem is significantly more tractable, \tou{} electricity pricing, charging slot capacity, and grid load are usually considered.
Surprisingly, it appears to be uncommon to incorporate non-linear battery behavior into charge scheduling models.

Since the field of \EVSP{} has made a lot of progress over the last twenty years, we highlight related papers featuring non-linear battery behavior or considering grid load, see \Cref{tab:literatureReview} for a summary.
To the best of our knowledge, no paper so far considers grid load and non-linear charging.
Grid load with linear charging is examined by \cite{MO19} and \cite{WLLJ22}.
\citeauthor{MO19} explicitly model parking spot assignment in the depot via a time discretization such that earlier arriving buses can potentially block later ones.
They track the total power draw along the time steps which enables them to consider \tou{} pricing and grid load.
The model is solved sequentially by heuristically finding a vehicle schedule, then the corresponding parking spot assignment and lastly the charging schedule.
\citeauthor{WLLJ22} consider a version without parking order and where buses are always charged fully, and present a column generation solving approach.

There are three approaches to handle non-linear battery behavior in the literature, which are energy state expansion, linear spline approximations, and algorithms that can consider arbitrary charge models exactly.
For a comparative overview see \Cref{tab:literatureReview}.

\begin{table}[htb]
    \centering
    \begin{tabular}{l | *{8}{c}}
        \toprule
        Reference & \makecell[c]{charge\\model} & \makecell[c]{not\\lin.} & \makecell[c]{part.\\chrg.} & \makecell[c]{dyn.\\rate} & \makecell[c]{\tou\\prices} & \makecell[c]{slot\\capa.} & \makecell[c]{grid\\load} & \makecell[c]{mix\\fleet} \\
        \midrule
        \cite{KAH17} & \rmexpansion & $\checkmark$ & $\checkmark$ &  & $\checkmark$ &  &  & \\
        \cite{MGMV17} & \rmspline & $\checkmark$ & $\checkmark$ &  &  &  &  & \\
        \cite{FMJL19} & \rmspline & $\checkmark$ & $\checkmark$ &  &  &  &  & \\
        \cite{LLX19} & \rmexpansion & $\checkmark$ &  &  &  & $\checkmark$ &  & \\
        \cite{MO19} & \rmlinear &  & $\checkmark$ & $\checkmark$ & $\checkmark$ & $\checkmark$ & $\checkmark$ & \\
        \cite{OK20} & \rmalgo & $\checkmark$ & $\checkmark$ &  &  &  &  & \\
        \cite{KGM21} & \rmspline & $\checkmark$ & $\checkmark$ &  &  & $\checkmark$ &  & \\
        \cite{Lee21} & \rmalgo & $\checkmark$ & $\checkmark$ &  &  &  &  & \\
        \cite{LDL21} & \rmspline & $\checkmark$ & $\checkmark$ &  &  &  &  & \\
        \cite{ZLTDCGL21} & \rmspline & $\checkmark$ & $\checkmark$ &  &  &  &  & \\
        \cite{ZWQ21} & \rmspline & $\checkmark$ &  &  &  & $\checkmark$ &  & \\
        \cite{FJML22} & \rmspline & $\checkmark$ & $\checkmark$ &  &  & $\checkmark$ &  & \\
        \cite{WLLJ22} & \rmlinear &  &  &  & $\checkmark$ & $\checkmark$ & $\checkmark$ & \\
        \cite{ZMO22} & \rmspline & $\checkmark$ & $\checkmark$ &  &  &  &  & \\
        \cite{DEG23} & \rmalgo & $\checkmark$ & $\checkmark$ &  &  &  &  & \\
        \cite{KC23} & \rmspline & $\checkmark$ & $\checkmark$ &  &  &  &  & \\
        \cite{VLD23} & \rmexpansion & $\checkmark$ & $\checkmark$ &  &  & $\checkmark$ &  & $\checkmark$ \\
        \midrule
        This Work & \rmrate & $\checkmark$ & $\checkmark$ & $\checkmark$ & $\checkmark$ & $\checkmark$ & $\checkmark$ & $\checkmark$ \\
        \bottomrule
    \end{tabular}
    \caption{Presented are publications and the features of the proposed models compared to \features{} specified in \Cref{sec:problemdef}. For the charge model column we write \rmlinear{} for linear charging, \rmexpansion{} for energy state expansion, \rmspline{} for the linear spline approach, \rmalgo{} (exact in algorithm) for considering a charge curve exactly within the solving algorithm, and \rmrate{} for our charge increment domain approach, see \Cref{sec:CID}.}
    \label{tab:literatureReview}
\end{table}

\subsection{Modeling Charging via Energy State Expansion}
\Citet{KAH17} propose energy state expansion in analogy to the well-known time expansion.
They consider a bus scheduling graph with a node for every combination of trip and discrete \soc.
Charging is fully encoded on arcs connecting trip nodes from lower to higher charge states.
Energy prices as well as battery fade are part of the arc costs.
The problem is solved using column generation in combination with Lagrangian relaxation.

\cite{LLX19} pursue a similar approach with a fully time-space-energy-expanded network, but consider different side constraints.
They determine the frequency at which trips should be repeated to satisfy a given passenger demand and the number of required chargers at the terminals.
They solve small instances with a commercial MILP solver.

To consider mixed electric fleets and a fixed number of charging slots, \cite{VLD23} use a time-and-energy-expanded network where charging slot occupation is tracked along discrete time steps.
They propose two rounding schemes for charge state propagation yielding primal and dual bounds.
They solve the problem with column generation and a diving heuristic to produce integer solutions.

\subsection{Modeling Charging via Linear Spline Approximation}
A series of papers by \cite{MGMV17}, \cite{FMJL19}, and \cite{FJML22} develop a local search heuristic based on a piecewise linear spline interpolation of a concave \emph{charge curve}, which can be incorporated into mixed-integer programs for the \EVSP{} via standard techniques.
The charge curve maps time spent charging an initially empty battery to the resulting \soc.
The local search produces vehicle course candidates over which a set partition problem is solved to obtain a final vehicle schedule. 
For each vehicle course, a cost optimal charge schedule has to be found to determine energy-feasibility and the total objective cost.
This \emph{fixed route vehicle charging problem} (\FRVCP) is NP-hard, so a labeling algorithm with a recharge event insertion heuristics is developed.
The shape of the concave piecewise linear charge curve yields efficient dominance criteria for the labels.

\cite{KGM21} further extend \cite{FMJL19} to logistic settings where operators may wish to charge at publicly accessible infrastructure owned by other parties, so that there is no control over the availability of charging slots.
They propose a benders-based branch-and-cut algorithm which uses a modified version of the \FRVCP{} labeling algorithm, accounting for potentially having to wait at charge stations by shifting the label function accordingly.

Approximating the charge curve by a linear spline has been adopted by a majority of subsequent papers on \EVSP{} with non-linear charging and thus can be considered the current state-of-the-art (see \Cref{tab:literatureReview}).
\cite{LDL21} develop a labeling algorithm for the pricing problem of a column generation approach using recursive functions derived from the linear spline.
\cite{ZLTDCGL21} limit the number of admissible line changes and present an adaptive large neighborhood search.
\cite{ZWQ21} develop a branch-and-price algorithm considering charger slot capacity and battery fade under a full recharge policy.
\cite{ZMO22} also consider battery fade and solve the resulting MILP including symmetry breaking inequalities with a commercial solver.
\cite{KC23} extend the model and algorithm from \cite{MGMV17} to also include non-linear discharging.

We discuss the linear spline modeling choice in \Cref{sec:rechargeModeling}.

\subsection{Considering Arbitrary Charge Models Exactly in Algorithms}
\cite{OK20} present a computational study to illustrate the importance of considering the non-linear battery behavior.
They use a greedy construction heuristic to obtain solutions with some backtracking to insert charging events.
This enables them to consider arbitrary recharge models exactly and they show that simple battery models, in particular those assuming a linear relationship between charge duration and \soc{}, can yield solutions that are infeasible in practice.

\cite{Lee21} proposes a branch-and-price algorithm for the \EVSP{} setting considered in \cite{MGMV17} and \cite{FMJL19} that can handle arbitrary charge curves as long as they are concave.
\citeauthor{Lee21} chooses an objective function that minimizes the total travel and charge time, so an optimal solution will only charge for as long as is necessary to cover the next sequence of customers.
As such, for every sequence of costumers between charging events a binary decision variable is introduced with the corresponding minimum required charge time as the objective coefficient.
The large number of variables then necessitates a branch-and-price framework and the charge curve is evaluated only upon checking the optimality criteria.
This approach seems incompatible with applications where an optimal solution is not guaranteed to minimize the total charge time, as is the case for the \EBSP{} as introduced in \Cref{sec:problemdef}.

\cite{DEG23} propose a branch-and-check algorithm for a factory in-plant electric tow train setting.
A master IP produces vehicle schedules and the subproblem explicitly computes the charge states under a policy where the tow trains are charged whenever and as much as possible.
The \soc{} values are computed exactly from arbitrary charge curves.
If the vehicle schedule turns out to not be energy-feasible, corresponding subtour elimination cuts are introduced into the main problem.

\section{Modeling the Recharge Process} \label{sec:rechargeModeling}
\begin{figure}[htb]
    \includegraphics[width=\textwidth]{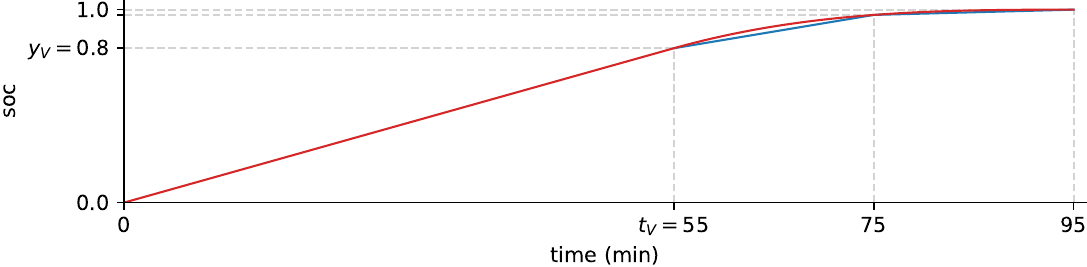}
    \caption{Example charge curve with commonly used linear spline interpolation.
    Based on real bus fast charging data.
    The interpolation has two segments for the CV phase.}
    \label{fig:ExampleChargeCurve}
\end{figure}

The literature describes the non-linear charging behavior of batteries by assuming the existence of a \emph{charge curve} $\chrgcurveMax: [0,\infty) \to [0,1]$ which maps time spent charging an initially empty battery to the resulting \soc.
Based on the \emph{Constant Current - Constant Voltage} (\cccv) charging scheme common for the lithium-ion battery type dominating the market \citep{PJLV17,OK20}, the charge curve has an initial linear segment until a point $\chrgcurveMax(\CCCVbreakTime) = \CCCVbreakSoC$, then it concavely and monotonically grows towards the maximum capacity, see \Cref{fig:ExampleChargeCurve} for an illustration.

\begin{definition}
    \label{def:classicCurve}
    We call $\chrgcurveMax: [0,\infty) \to [0,1]$ a \emph{charge curve} if
    \begin{enumerate}[(i)]
        \item \label{def:classicCurve:curveMonotone} there is an interval $[0, \chrgtimeMax]$ on which $\chrgcurveMax$ is strictly monotonically increasing,
        \item \label{def:classicCurve:curveBijective} $\chrgcurveMax(\timevar) = 1$ holds for $\timevar \geq \chrgtimeMax$,
        \item \label{def:classicCurve:derivative} and the derivative $\chrgcurveMax'$ exists almost everywhere and is monotonically non-increasing.
    \end{enumerate}
\end{definition}

From a charge curve one usually has to determine the time required to charge from an initial to some final \soc, or to find the amount of replenished driving range upon charging for a given duration from some initial charge state.
This can be formalized as follows.

\begin{definition}
    \label{def:operators}
    Let $\chrgcurveMax$ be a charge curve.
    Then an inverse $\chrgcurveMax^{-1}: [0,1] \to [0, \chrgtimeMax]$ exists and the \emph{charge duration operator} $\timeOp{}$ yields a \emph{charge duration function} $\timeOp{\chrgcurveMax}$ defined by
    \begin{equation}
        \timeOp{\chrgcurveMax}: \setEnc{(\initSoC, \finalSoC) \in [0,1]^2 \mid \initSoC \leq \finalSoC} \to [0,\infty),~ \timeOp{\chrgcurveMax}(\initSoC, \finalSoC) = \chrgcurveMax^{-1}(\finalSoC) - \chrgcurveMax^{-1}(\initSoC)
    \end{equation}
    and the \emph{charge increment operator} $\diffOp{}$ yields a \emph{charge increment function} $\diffOp{\chrgcurveMax}$ defined by
    \begin{equation}
        \diffOp{\chrgcurveMax}: [0,1] \times [0,\infty) \to [0,1],~ \diffOp{\chrgcurveMax}(\chrgvar, \timevar) = \chrgcurveMax(\chrgcurveMax^{-1}(\chrgvar) + \timevar) - \chrgvar\text{.}
    \end{equation}
\end{definition}

\begin{lemma}
    \label{lemma:diffOpMonotone}
    $\diffOp{\chrgcurveMax}(\chrgvar, \timevar)$ is monotonically non-increasing in $\chrgvar$.
\end{lemma}
\proof
\begin{equation}
    \frac{\partial}{\partial \chrgvar} \diffOp{\chrgcurveMax}(\chrgvar, \timevar) = \frac{\chrgcurveMax'(\chrgcurveMax^{-1}(\chrgvar) + \timevar)}{\chrgcurveMax'(\chrgcurveMax^{-1}(\chrgvar))} - 1 \leq 0
\end{equation}
by \Cref{def:classicCurve}.\ref{def:classicCurve:derivative} and $\timevar \geq 0$. \myendproof

We will focus on the increment operator in this paper, one can argue analogously for the duration operator.
The point of either operator is to obtain certificates asserting that classically feasible vehicle schedules, or in general any route or path an electric vehicle might take, together with some charge schedule are also \emph{energy-feasible}, a notion we have to formally define.
Consider any vehicle course $\course$ and recall that it is an alternating sequence of $m$ blocks and $m - 1$ recharge events, assuming the bus can reach the depot from its last trip without having to charge again.
For every block $\block_i$ on $\course$, enumerate the start and end terminals of each trip such that we can write
\begin{equation}
    \block_i = (\trip_{k_{i-1} + 1}^s, \trip_{k_{i-1} + 1}^e, \trip_{k_{i-1} + 2}^s, \trip_{k_{i-1} + 2}^e, \dots, \trip_{k_i}^s, \trip_{k_i}^e)
\end{equation}
with $k_0 = 0$ and $i\in [m] \defeq \setEnc{1,\dots,m}$.
Then together with the initial departure at the depot $\dep^s$, the arrival and departure for each recharge event $(\chrgslot_i^s, \chrgslot_i^e)$, $i \in [m - 1]$, and the final arrival at the depot $\dep^e$, we can express the vehicle course as a tuple
\begin{multline}
    \course = (\dep^s) \times \block_1 \times (\chrgslot_1^s, \chrgslot_1^e) \times \block_2 \times (\chrgslot_2^s, \chrgslot_2^e) \times \dots \times (\chrgslot_{m-1}^s, \chrgslot_{m-1}^e) \times \block_m \times (\dep^e) \\
    = (\dep^s, \trip_1^s, \trip_1^e, \dots, \trip_{k_1}^s, \trip_{k_1}^e, \chrgslot_1^s, \chrgslot_1^e, \trip_{k_1 + 1}^s, \trip_{k_1 + 1}^e \dots, \trip_{k_m}^s, \trip_{k_m}^e, \dep^e)\text{.}
\end{multline}
Furthermore, a charge time window of duration $\timevar_i$ is associated to each recharge event $(\chrgslot_i^s, \chrgslot_i^e)$.
We can then propagate charge states along the vehicle course $\course$ using an iterative map
\begin{equation}
    \label{eq:chrgPropagation}
    \chrgvar: \course \to \mathbb{R},\quad \begin{cases}
        \dep^s & \mapsto 1 \\
        \chrgslot_i^e & \mapsto \chrgvar(\chrgslot_i^s) + \diffOp{\chrgcurveMax}(\chrgvar(\chrgslot_i^s), \timevar_i) \\
        \course_j & \mapsto \chrgvar(\course_{j - 1}) - \consumption(\course_{j - 1}, \course_j)
    \end{cases}
\end{equation}
where $\consumption(.,.)$ is a non-negative function giving the relative energy consumption scaled to $[0,1]$ between two locations.
Note that arrival and departure locations for depots and charge stations may differ and that $\consumption(.,.)$ does not need to be commutative, but it should satisfy the triangle inequality $\consumption(a,b) \leq \consumption(a,c) + \consumption(c,b)$.

While the laws of physics dictate that the \soc{} has to remain non-negative, operators usually want to maintain a minimum energy level for flexible vehicle disposition.
This minimum charge state should at least allow a bus to reach the nearest depot or charge station from its current planned location.
\begin{definition}
    Let $\course$ be a vehicle course and let $\minReqChrg(\course_j) \geq 0$ be the minimum required \soc{} to reach the nearest depot or charge station from course location $\course_j$.
    We say $\course$ is \emph{energy-feasible} (under $\chrgcurveMax$) if and only if $\chrgvar(\course_j) \geq \minReqChrg(\course_j)$ for all $j\in [\abs{\course}]$.
    A vehicle schedule is energy-feasible if all of its vehicle courses are.
\end{definition}

Replacing the charge curve by any approximation introduces an error which propagates to the charge duration respectively charge increment operator and further onto the charge state map defined by \eqref{eq:chrgPropagation}.
As such, let $\approxOp{\chrgcurveMax}$ be a charge curve which in some sense approximates the exact curve $\chrgcurveMax$.
Note that a linear spline interpolation of $\chrgcurveMax$ complies with this.
Then let
\begin{equation}
    \errorFunc(\chrgvar, \timevar) \defeq \diffOp{\approxOp{\chrgcurveMax}}(\chrgvar, \timevar) - \diffOp{\chrgcurveMax}(\chrgvar, \timevar)
\end{equation}
give the approximation error for the charge increment operator.
Replacing $\chrgcurveMax$ by $\approxOp{\chrgcurveMax}$ in \eqref{eq:chrgPropagation} yields an approximate charge propagation map $\approxOp{\chrgvar}$ and we obtain the stepwise approximation error along the course as
\begin{equation}
    \errorFunc(\course_j) \defeq \approxOp{\chrgvar}(\course_j) - \chrgvar(\course_j)\text{.}
\end{equation}
Consequently, a vehicle course is energy-feasible under $\chrgcurveMax$ if and only if
\begin{equation}
    \label{eq:approxFeasCrit}
    \approxOp{\chrgvar}(\course_j) \geq \minReqChrg(\course_j) + \errorFunc(\course_j) \text{ for all } j\in[\abs{\course}]
\end{equation}
holds.
While $\minReqChrg$ is known a priori, $\errorFunc$ depends on the vehicle course by
\begin{equation}
    \label{eq:errorPropagation}
    \errorFunc(\chrgslot_i^e) = \errorFunc(\chrgslot_{i-1}^e) + \diffOp{\approxOp{\chrgcurveMax}}(\approxOp{\chrgvar}(\chrgslot_i^s), \timevar_i) - \diffOp{\chrgcurveMax}(\chrgvar(\chrgslot_i^s), \timevar_i) = \approxOp{\chrgcurveMax}(\approxOp{\chrgcurveMax}^{-1}(\approxOp{\chrgvar}(\chrgslot_i^s)) + \timevar_i) - \chrgcurveMax(\chrgcurveMax^{-1}(\chrgvar(\chrgslot_i^s)) + \timevar_i)
\end{equation}
from error propagation at recharge events, which in general is difficult to capture.
Nevertheless, we have to study the numerical stability of $\approxOp{\chrgvar}$.
In particular, we have to be interested in the worst-case impact of an approximation error, its bounds, and how the state-of-the-art linear spline approximation of a charge curve behaves.
We will show that it is preferable to interpolate the charge increment function as described in \Cref{sec:CID} instead of the charge curve. 

\subsection{A Worst Case Instance}
We demonstrate how to construct a worst case instance for charge curve approximation in general.
\begin{definition}
    $\course$ is \emph{weakly energy-feasible} if
    $\approxOp{\chrgvar}(\course_j) \geq \minReqChrg(\course_j)$ holds for all $j\in[\abs{\course}]$.
\end{definition}
\begin{definition}
    $\course$ is \emph{strongly energy-feasible} if it is energy-feasible and also weakly so.
\end{definition}

\begin{theorem}\label{theorem:worstcase}
    Suppose assigning a unique bus to each trip yields an energy-feasible solution, i.e., for every trip $\trip \in \trips$ there is a depot $d$ such that $\consumption(d,\trip^s) + \consumption(\trip^s,\trip^e) + \consumption(\trip^e,d) \leq 1$.
    Then the minimum number of courses for a strongly energy-feasible vehicle schedule under $\chrgcurveMax$ and $\approxOp{\chrgcurveMax}$ is at most $\abs{\trips}$ times the minimum number of courses for an energy-feasible schedule and the bound is tight.
\end{theorem}
This is clearly the worst possible bound, but we will see that a corresponding worst-case instance has to have highly unrealistic properties.
\proof The inequality trivially holds by assumption. To show tightness we construct an instance where the optimal energy-feasible vehicle schedule consists of a single vehicle course, and the optimal strongly energy-feasible vehicle schedule requires $\abs{\trips}$ courses, one per trip.
Let $\course$ be such an energy-feasible course.
It is not weakly energy-feasible if there is some trip such that $\approxOp{\chrgvar}(\trip_j^e) < \minReqChrg(\trip_j^e)$ $\leq \chrgvar(\trip_j^e)$ holds.
Assuming we can not extend any earlier recharge events to compensate for the approximation error, we need a separate bus to service trip $\trip_j$.

\begin{lemma}\label{lemma:newbus}
    If energy consumption satisfies the triangle inequality, we can always send a replacement bus which will arrive with a \soc{} of at least $\chrgvar(\trip_j^s)$ at $\trip_j^s$.
\end{lemma}
\proof Analogously to how $\minReqChrg(\blankArg)$ denotes the lowest admissible \soc{} at a location, let $\maxChrg(\blankArg)$ denote the highest possible \soc{} a bus coming directly from any depot or charge station can have at the given location.
Further, let $\ell$ denote the depot or charge station preceding the vehicle block containing $\trip_j$.
If every possible replacement bus would arrive with a \soc{} strictly less than $\chrgvar(\trip_j^s)$, we get
\begin{equation}
    \maxChrg(\trip_j^s) < \chrgvar(\trip_j^s) =  \chrgvar(\ell) - \consumption(\ell, \trip_i^s) - \dots - \consumption(\trip_{j-1}^e, \trip_j^s) \leq 1 - \consumption(\ell, \trip_j^s) \leq \maxChrg(\trip_j^s) \text{,}
\end{equation}
a contradiction.\myendclaimproof

(\emph{Continuation of proof for \Cref{theorem:worstcase}})
If we can always send a replacement bus that has a sufficient \soc{} buffer at $\trip_j^s$, it will be able to service the rest of the course and the worst case bound on the fleet size increase would be $2$.
If, however, $\maxChrg(\trip_j^s) \approx \chrgvar(\trip_j^s)$ and there is another critical trip violating weak energy-feasibility, we might need a third replacement bus, and so forth.
The optimal energy-feasible course $\course$ of a worst-case instance therefore has to be an alternating sequence of trips and recharge events.

Moreover, to achieve the worst case bound of $\abs{\trips}$ courses, all already deployed buses must not be able to chain any of the recharge events on the original course $\course$.
For example, while we assume that the initial bus can not acquire enough charge under the approximation to service the second block, if it restores enough driving range to reach the second recharge event, it could replenish enough charge there to entirely finish the course $\course$ afterwards, giving us again a bound of $2$.
This has a few implications:
Every recharge event must take place at a unique location which is available exactly during the tight time window of the original course $\course$, otherwise a bus could simply wait for the next window to open and acquire enough charge to potentially resume the original course later.
This is theoretically achievable via feature \ref{feat:gridload} by setting the grid load limit to zero outside the required charge windows.
We do have real-life \EBSP{} instances where charging is at least not allowed in the evening.
The recharge locations then have to be spaced such that no deployed bus is able to charge enough under the approximation by hopping recharge events, so $\approxOp{\chrgvar}(\chrgslot_i^e) - \consumption(\chrgslot_i^e, \chrgslot_j^s) < \chrgvar(\chrgslot_j^s)$ for every later recharge window.
In particular,
\begin{multline}\label{eq:minStationDistance}
    \consumption(\chrgslot_{i-1}^e, \trip_i^s) + \consumption(\trip_i^s, \trip_i^e) + \consumption(\trip_i^e, \chrgslot_i^s) \geq \consumption(\chrgslot_{i-1}^e, \chrgslot_i^s) > \approxOp{\chrgvar}(\chrgslot_{i-1}^e) - \chrgvar(\chrgslot_i^s) \\
    = \errorFunc(\chrgslot_{i-1}^e) + \consumption(\chrgslot_{i-1}^e, \trip_i^s) + \consumption(\trip_i^s, \trip_i^e) + \consumption(\trip_i^e, \chrgslot_i^s)\text{,}
\end{multline}
that is, each trip has to approximately lie on the path between its adjacent charge stations.

\begin{figure}[htb]
    \resizebox{\textwidth}{!}{\begin{tikzpicture}[
    node/.style = {circle, draw=black, fill=white, inner sep = 0pt, minimum size=0.71cm}, 
    depot/.style = {regular polygon, regular polygon sides=4, draw=black, fill=white, inner sep = 0pt, minimum size=1cm},
    station/.style = {regular polygon, regular polygon sides=3, shape border rotate=180, draw=black, fill=white, inner sep = 0pt, minimum size=1cm},
    edge/.style = {decoration={markings,mark=at position 1 with {\arrow[scale=2,>=stealth]{>}}},postaction={decorate}}
    ]

    \begin{scope} [every node/.style={depot}]
        \node (d0) at (1,-1.732) {$\dep_1$};
        \node (d1) at (7,-1.732) {$\dep_2$};
        \node (d0again) at (15,-1.732) {$\dep_1$};
    \end{scope}
    \begin{scope} [every node/.style={station}]
        \node (s0) at (4,0) {$\chrgslot_1$};
        \node (s1) at (10,0) {$\chrgslot_2$};
    \end{scope}
    \begin{scope} [every node/.style={node}]
        \node (t0s) at (0,0) {$\trip_1^s$};
        \node (t0e) at (2,0) {$\trip_1^e$};

        \node (t1s) at (6,0) {$\trip_2^s$};
        \node (t1e) at (8,0) {$\trip_2^e$};

        \node (tns) at (14,0) {$\trip_n^s$};
        \node (tne) at (16,0) {$\trip_n^e$};
    \end{scope}

    \node (dots) at (12,0) {\dots};

    \draw[edge] (d0) to node[sloped, anchor=center, below] {$\frac{1}{3}$} (t0s);
    \draw[edge] (t0s) to node[sloped, anchor=center, above] {$\frac{1}{3} + \delta$} (t0e);
    \draw[edge] (t0e) to node[sloped, anchor=center, below] {$\frac{1}{3}$} (d0);

    \draw[edge] (t0e) to node[sloped, anchor=center, above] {$\frac{1}{3} - \delta$} (s0);
    \draw[edge] (s0) to node[sloped, anchor=center, above] {$\frac{1}{3} + \delta$} (t1s);

    \draw[edge] (d1) to node[sloped, anchor=center, below] {$\frac{1}{3}$} (t1s);
    \draw[edge] (t1s) to node[sloped, anchor=center, above] {$\frac{1}{3}$} (t1e);
    \draw[edge] (t1e) to node[sloped, anchor=center, below] {$\frac{1}{3}$} (d1);

    \draw[edge] (t1e) to node[sloped, anchor=center, above] {$\frac{1}{3} - \delta$} (s1);
    \draw[edge] (s1) to node[sloped, anchor=center, above] {$\frac{1}{3} + \delta$} (dots);
    \draw[edge] (dots) to node[sloped, anchor=center, above] {$\frac{1}{3} + \delta$} (tns);

    \draw[edge] (d0again) to node[sloped, anchor=center, below] {$\frac{1}{3}$} (tns);
    \draw[edge] (tns) to node[sloped, anchor=center, above] {$\frac{1}{3}$} (tne);
    \draw[edge] (tne) to node[sloped, anchor=center, below] {$\frac{1}{3}$} (d0again);
\end{tikzpicture}}
    \caption{Spatial network of a worst case \EBSP{} instance for charge curve approximation with relative energy consumption and $n$ trips, $n-1$ recharge stations, and $n-1$ depots. 
    Trips $\trip_1$ and $\trip_n$ are both attached to $\dep_1$ since buses have to return to their depot of origin.}
    \label{fig:underestWorstCase}
\end{figure}
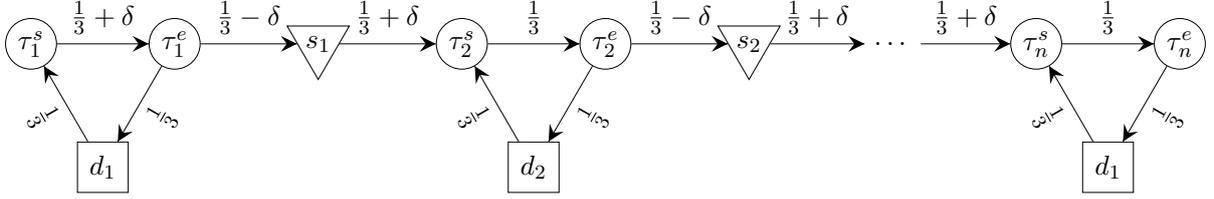

A worst-case instance therefore has to be structured as in \Cref{fig:underestWorstCase} where we take any $0 > \delta > \errorFunc$ for an approximation error $\errorFunc$ charging an empty battery to get (almost) full.

Note that since $\trip_n$ has to be attached to $\dep_1$, the network has to have some curvature, which in turn is limited by \eqref{eq:minStationDistance} so that no considerable shortcuts emerge.
The instance therefore has to have a minimum number of trips which depend on $\errorFunc$.

If each recharge event allows a bus to exactly fully recharge its battery, then a single bus can service every trip in order.
However, under any (underestimating) approximation error, we have $\approxOp{\chrgvar}(\chrgslot_1^e) = 1 + \epsilon < 1 + \delta$, which is not enough to drive the second trip or to reach any other recharge event.
A second bus can be deployed from $\dep_2$ or $\dep_1$.
Coming from $\dep_2$, it can not reach $\chrgslot_2$ while servicing $\trip_2$.
From $\dep_1$, it could take the first charge window and achieve a full battery since its arrival charge state is higher than if it serviced $\trip_1$.
But then it arrives empty at $\chrgslot_2$, so we need a third bus for $\trip_3$, and so forth.\myendproof

\begin{corollary}\label{corollary:worstcase}
    The minimum number of courses for an energy-feasible vehicle schedule is at most $\abs{\trips}$ times the minimum number of courses for a weakly energy-feasible schedule and the bound is tight.
\end{corollary}
\proof The same worst-case instance works, except we exchange the roles of $\chrgvar$ and $\approxOp{\chrgvar}$, i.e., for a critical trip we have $\chrgvar(\trip_j^e) < \minReqChrg(\trip_j^e) \leq \approxOp{\chrgvar}(\trip_j^e)$ and thus $0 < \delta < \epsilon$, so the signs of $\delta$ in the energy consumption in the network in  \Cref{fig:underestWorstCase} have to be flipped.
The charge windows are now exactly such that the approximation overestimation allows a single bus to service all trips, but under the exact curve the windows are too tight and we need $\abs{\trips}$ buses.\myendproof

We can reasonably expect no real-life \EBSP{} instance to exhibit this worst-case behavior, however, it is certainly possible that an approximate charge curve could have an impact on the minimum fleet size.
Any underestimation can lead to the loss of optimal courses, whereas any overestimation can produce a false-positive energy-feasibility certificate.

\subsection{Bounds on the Error Propagation}
A straightforward consequence of \Cref{theorem:worstcase} and \Cref{corollary:worstcase} is that we should prefer charge curve approximations such that $\approxOp{\chrgvar} \leq \chrgvar$ holds along any possible vehicle course, i.e., that weak energy-feasibility implies strong energy-feasibility.
We can then at least avoid \Cref{corollary:worstcase} and be certain that we never obtain a false-positive energy-feasibility certificate for a vehicle course.
Fortunately, requiring $\diffOp{\approxOp{\chrgcurveMax}}(\chrgvar, \timevar) \leq \diffOp{\chrgcurveMax}(\chrgvar, \timevar)$ suffices:

For the following, we write $\course_i \preceq \course_j$ to indicate that vehicle course element $\course_i$ happens before $\course_j$ or that $\course_i = \course_j$.
\begin{proposition}
    \label{prop:propagationErrorUB}
    If $\errorFunc(\chrgvar, \timevar) = \diffOp{\approxOp{\chrgcurveMax}}(\chrgvar, \timevar) - \diffOp{\chrgcurveMax}(\chrgvar, \timevar) \leq 0$ for all $\chrgvar \in [0,1]$ and $\timevar \in [0, \chrgtimeMax]$, then for every element $\course_j \preceq \chrgslot_1^s$ before the first recharge event of a vehicle course $\course$ we have $\errorFunc(\course_j) = 0$, and for the remaining elements $\chrgslot_1^e \preceq \course_j$ we have $0 \geq \errorFunc(\chrgvar(\chrgslot_{\numChrgEvents(j)}^s), \timevar_{\numChrgEvents(j)}) \geq \errorFunc(\course_j)$ where $\numChrgEvents(j) \defeq \abs{\setEnc{i \in [m-1] \mid \chrgslot_i^e \preceq \course_j}}$ gives the number of recharge events before $\course_j$ and thus the index of the event immediately preceding the block containing $\course_j$.
\end{proposition}
\proof Since $\errorFunc$ remains constant on individual blocks, it suffices to examine the recharge events, so the proof is by induction on their number in the course $\course$.
By definition, $\approxOp{\chrgvar}(\course_j) = \chrgvar(\course_j)$ on $(\dep^s) \times \block_1 \times (\chrgslot_1^s)$ and $\errorFunc(\chrgslot_1^e) = \diffOp{\approxOp{\chrgcurveMax}}(\chrgvar(\chrgslot_1^s), \timevar_1) - \diffOp{\chrgcurveMax}(\chrgvar(\chrgslot_1^s), \timevar_1) = \errorFunc(\chrgvar(\chrgslot_1^s), \timevar_1) \leq 0$, which is propagated onto $(\chrgslot_1^e) \times \block_2 \times (\chrgslot_2^s)$.

Now consider $\chrgslot_i^e$ with $i \geq 2$.
By induction, $0 \geq \errorFunc(\chrgslot_{i-1}^e) = \errorFunc(\chrgslot_i^s)$, thus $\approxOp{\chrgvar}(\chrgslot_i^s) \leq \chrgvar(\chrgslot_i^s)$ and then
\begin{multline}
    \errorFunc(\chrgslot_i^e) = \approxOp{\chrgcurveMax}(\approxOp{\chrgcurveMax}^{-1}(\approxOp{\chrgvar}(\chrgslot_i^s)) + \timevar_i) - \chrgcurveMax(\chrgcurveMax^{-1}(\chrgvar(\chrgslot_i^s)) + \timevar_i) \leq \approxOp{\chrgcurveMax}(\approxOp{\chrgcurveMax}^{-1}(\chrgvar(\chrgslot_i^s)) + \timevar_i) - \chrgcurveMax(\chrgcurveMax^{-1}(\chrgvar(\chrgslot_i^s)) + \timevar_i) \\
    = \diffOp{\approxOp{\chrgcurveMax}}(\chrgvar(\chrgslot_i^s), \timevar_i) + \chrgvar(\chrgslot_i^s) - \diffOp{\chrgcurveMax}(\chrgvar(\chrgslot_i^s), \timevar_i) - \chrgvar(\chrgslot_i^s) = \errorFunc(\chrgvar(\chrgslot_i^s), \timevar_i)
\end{multline}
by \eqref{eq:errorPropagation} and monotonicity of $\approxOp{\chrgcurveMax}$.\myendproof

\begin{corollary}
    If $\errorFunc(\chrgvar, \timevar) \geq 0$ for all $\chrgvar \in [0,1]$ and $\timevar \in [0, \chrgtimeMax]$, then $\errorFunc(\course_j) = 0$ for every $\course_j \preceq \chrgslot_1^s$, and $0 \leq \errorFunc(\chrgvar(\chrgslot_{\numChrgEvents(j)}^s), \timevar_{\numChrgEvents(j)}) \leq \errorFunc(\course_j)$ for the remaining elements $\chrgslot_1^e \preceq \course_j$.
\end{corollary}

\begin{proposition}
    If $\errorFunc(\chrgvar, \timevar) \leq 0$ for all $\chrgvar \in [0,1]$ and $\timevar \in [0, \chrgtimeMax]$, then for every element $\course_j$ of a vehicle course $\course$,
    \begin{equation}
        \errorFunc(\course_j) \geq \sum_{\chrgslot_i^e \in \course:~\chrgslot_i^e \preceq \course_j} \errorFunc(\chrgvar(\chrgslot_i^s), \timevar_i)\text{.}
    \end{equation}
\end{proposition}
\proof We have already seen that $\errorFunc(\course_j) = 0$ on the first block $(\dep^s) \times \block_1 \times (\chrgslot_1^s)$ and on the second block $(\chrgslot_1^e) \times \block_2 \times (\chrgslot_2^s)$ we have $\errorFunc(\course_j)$ $= \errorFunc(\chrgvar(\chrgslot_1^s), \timevar_1)$.
Consider $i \geq 2$ and note that by \Cref{prop:propagationErrorUB}, $\approxOp{\chrgvar}(\chrgslot_i^s) \leq \chrgvar(\chrgslot_i^s)$.
Further recall \Cref{lemma:diffOpMonotone} that $\diffOp{\approxOp{\chrgcurveMax}}$ is monotonically non-increasing.
Induction on $i$ then yields
\begin{multline}
    \errorFunc(\chrgslot_i^e) = \errorFunc(\chrgslot_{i-1}^e) + \diffOp{\approxOp{\chrgcurveMax}}(\approxOp{\chrgvar}(\chrgslot_i^s), \timevar_i) - \diffOp{\chrgcurveMax}(\chrgvar(\chrgslot_i^s), \timevar_i) \\
    \geq \sum_{j = 1}^{i - 1} \errorFunc(\chrgvar(\chrgslot_j^s), \timevar_j) + \diffOp{\approxOp{\chrgcurveMax}}(\chrgvar(\chrgslot_i^s), \timevar_i) - \diffOp{\chrgcurveMax}(\chrgvar(\chrgslot_i^s), \timevar_i) = \sum_{j = 1}^{i} \errorFunc(\chrgvar(\chrgslot_j^s), \timevar_j)
\end{multline}
which proves the claim as $\errorFunc$ is constant along blocks.\myendproof

\begin{corollary}
    If $\errorFunc(\chrgvar, \timevar) \geq 0$ for all $\chrgvar \in [0,1]$ and $\timevar \in [0, \chrgtimeMax]$, then for every element $\course_j$ of a vehicle course $\course$,
    \begin{equation}
        \errorFunc(\course_j) \leq \sum_{\chrgslot_i^e \in \course:~\chrgslot_i^e \preceq \course_j} \errorFunc(\chrgvar(\chrgslot_i^s), \timevar_i)\text{.}
    \end{equation}
\end{corollary}

\begin{corollary}
    \label{cor:totalErrorAlongCourse}
    If $\errorFunc(\chrgvar, \timevar)$ has the same sign for all arguments, then $\errorFunc(\course_j)$ shares that sign along the entire course $\course$ and $\abs{\errorFunc(\course_j)} \leq \numChrgEvents(j) \, \lVert \diffOp{\approxOp{\chrgcurveMax}} - \diffOp{\chrgcurveMax} \rVert$ where $\norm{.}$ is the maximum norm.
\end{corollary}
In other words, if $\diffOp{\approxOp{\chrgcurveMax}}$ always underestimates (overestimates) $\diffOp{\chrgcurveMax}$, then so does $\approxOp{\chrgvar}(\course_j)$ underestimate (overestimate) $\chrgvar(\course_j)$.
\begin{corollary}
    Weak energy-feasibility implies strong energy-feasibility if $\diffOp{\approxOp{\chrgcurveMax}} \leq \diffOp{\chrgcurveMax}$.
\end{corollary}

\subsection{Linear Spline Charge Curve Approximation}
From the previous discussion it is clear that any charge curve approximation for which neither $\diffOp{\approxOp{\chrgcurveMax}} \leq \diffOp{\chrgcurveMax}$ nor $\diffOp{\approxOp{\chrgcurveMax}} \geq \diffOp{\chrgcurveMax}$ holds is unreliable, as it can arbitrarily over- and underestimate charge states along vehicle courses.
Unfortunately, linear spline approximations of the charge curve are such unreliable approximations.

\begin{figure}[htb]
    \includegraphics[width=0.49\textwidth]{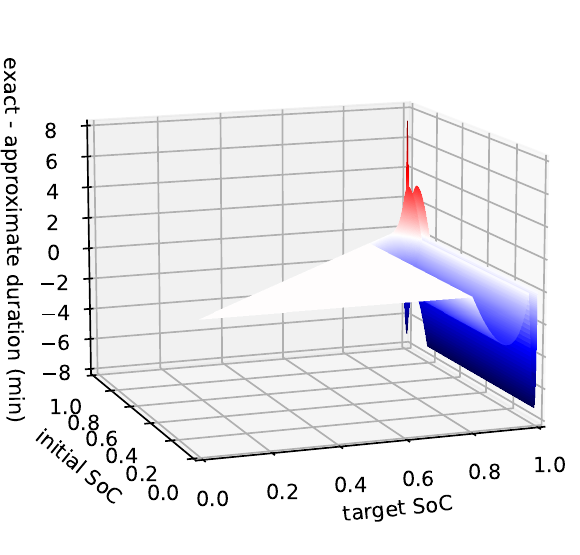}
    \includegraphics[width=0.49\textwidth]{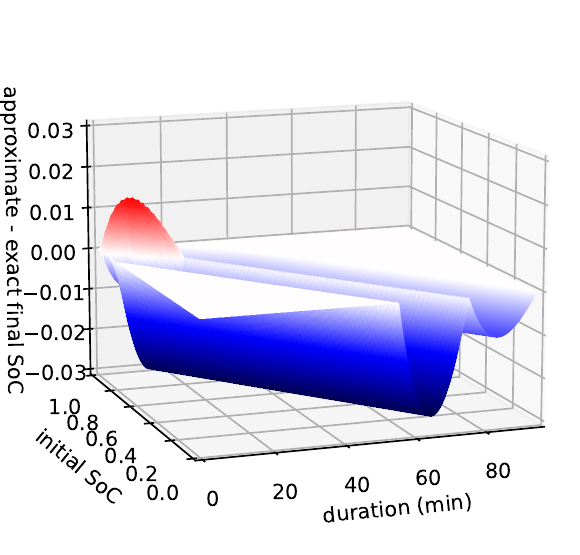}
    \caption{$\timeOp{\chrgcurveMax} - \timeOp{\linSpline{\chrgcurveMax}}$ on the left, $\diffOp{\linSpline{\chrgcurveMax}} - \diffOp{\chrgcurveMax}$ on the right for the example charge curve and linear spline approximation in \Cref{fig:ExampleChargeCurve}.}
    \label{fig:ApproxError3D}
\end{figure}

\begin{definition}
    Given a function $f: [a,b] \to \mathbb{R}$ and some grid on $[a,b]$, let $\linSpline{f}$ denote the (unique) piecewise linear spline interpolation of $f$ on that grid.
\end{definition}
\begin{proposition}\label{linSplineOscillates}
    Given a linear spline interpolation $\approxOp{\chrgcurveMax} = \linSpline{\chrgcurveMax}$ of a charge curve $\chrgcurveMax$ on $[0, \chrgtimeMax]$, there exist arguments $\chrgvar$ and $\timevar$ for which $\errorFunc(\chrgvar, \timevar) < 0$ and and there exist arguments for which $\errorFunc(\chrgvar, \timevar) > 0$.
\end{proposition}
\proof \cite{LBWor24}. \myendproof
So if a linear spline interpolation $\linSpline{\chrgcurveMax}$ of $\chrgcurveMax$ is used as suggested by \cite{MGMV17}, $\diffOp{\linSpline{\chrgcurveMax}}$ (and $\timeOp{\linSpline{\chrgcurveMax}}$) will \emph{both over- and underestimate} their exact counterparts, see \Cref{fig:ApproxError3D} for an example visualization of this behavior.

To provide an intuitive understanding why this happens, observe that the derivative of the charge curve states the applied charge rate at time $\timevar$.
Even though we can guarantee $\linSpline{\chrgcurveMax} \leq \chrgcurveMax$ for the linear spline interpolation of a concave function, the derivative of a piecewise linear function is piecewise constant, and hence $(\linSpline{\chrgcurveMax})'$ periodically overestimates the exact charge rate given by $\chrgcurveMax'$, such that we can find input values for which the average charge rate is overestimated.
This behavior has been observed by \cite{ZLTDCGL21} in their computational study and it even occurs when using a charge curve underestimator that is just the linear function through $(0,0)$ and $(\chrgtimeMax,1)$.
We formalized this result as an inherent property of piecewise linear charge curves in \citep{LBWor24} because the literature at large seems unaware of \Cref{linSplineOscillates} and appears to operate under the incorrect assumption that a spline underestimation of the charge curve results in an underestimation of the charge increment function, i.e.,  
\begin{equation}\label{eq:incorrectAssumption}
    \linSpline{\chrgcurveMax}(\timevar) \leq \chrgcurveMax(\timevar) ~\forall\, \timevar \in [0, \chrgtimeMax] \Longrightarrow \errorFunc(\chrgvar, \timevar) = \diffOp{\linSpline{\chrgcurveMax}}(\chrgvar, \timevar) - \diffOp{\chrgcurveMax}(\chrgvar, \timevar) \leq 0 ~\forall\, \chrgvar \in [0,1], \timevar \in [0, \chrgtimeMax]\text{.}
\end{equation}

\section{Linearizing Dynamic Charging via the Charge Increment Domain}\label{sec:CID}
We recapitulate and extend upon \citep{LBWatmos23} in this section.

The previous discussion and the emerging requirements of public transport operators outlined in \Cref{sec:intro} suggest to develop recharge models that focus on the charge rate.

Common battery types have a maximum safe charging rate that decreases with the present \soc, which necessitates the \cccv{} charging scheme \citep{PJLV17}.
In other words, there is a maximum admissible charge rate depending solely on how full the battery already is.
\begin{definition}\label{def:chargingProfile}
    A \emph{charging power profile} is a function $\ODEfunc: [0,1] \to [0,1]$ that maps the (relative) battery \soc{} to the (relative) maximal charge rate.
    There is some break point $\CCCVbreakSoC \in (0,1)$ such that $\ODEfunc$ is of the form
    \begin{equation}
        \ODEfunc(\chrgvar) = \begin{cases}
            \ODEfuncCC, & \chrgvar < \CCCVbreakSoC \\
            \ODEfuncCV(\chrgvar), & \chrgvar \geq \CCCVbreakSoC,
        \end{cases}
    \end{equation}
    where $\ODEfuncCV$ is differentiable, monotonically non-increasing, and satisfies $\ODEfuncCV(\CCCVbreakSoC) = \ODEfuncCC$ as well as $\ODEfuncCV(1) = 0$.
\end{definition}

\begin{figure}[htb]
    \includegraphics[width=\textwidth]{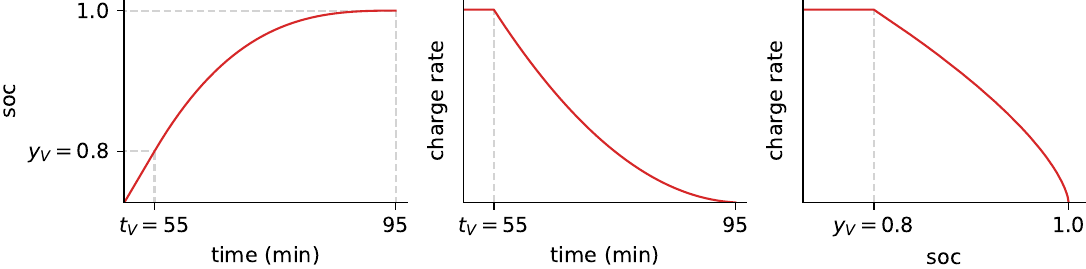}
    \caption{From left to right: CV segment of the charge curve from \Cref{fig:ExampleChargeCurve}, its derivative, and its power profile.}
    \label{fig:chargeCurveDeriv}
\end{figure}

A charger does not have to apply the highest possible charge rate.
In principle, the difference between slow and fast charging is the initial constant current.
High currents reach the \cccv{} breakpoint sooner, age the battery faster, and may contribute to grid load peaks.
Considering increasing adoption of active charge management tools which can dynamically throttle the charging current, we want a charge curve to be a description of an arbitrary recharge process and so we generalize \Cref{def:classicCurve}.
\begin{definition}\label{def:chargeCurve}
    Given a charging power profile $\ODEfunc$, a \emph{charge curve} describes an admissible charging process mapping time spent charging an initially empty battery to the resulting \soc{}.
    It is a differentiable function $\chrgcurve: [0,\infty) \to [0,1]$ satisfying
    \begin{enumerate}[(i)]
        \item \label{enum:ODEleft} $\chrgcurve(0) = 0$,
        \item \label{enum:ODEright} there exists $\chrgtimeMax > 0$ such that $\chrgcurve(\timevar) = 1$ for $\timevar \geq \chrgtimeMax$,
        \item \label{enum:ODE} and $0 < \chrgcurve'(\timevar) \leq \ODEfunc(\chrgcurve(\timevar))$ for all $\timevar \in [0, \chrgtimeMax)$.
    \end{enumerate}
    Let $\curves(\ODEfunc)$ denote the set of all charge curves derived from $\ODEfunc$.
\end{definition}
By relaxing $0 < \chrgcurve'(\timevar)$ we can account for temporarily suspending charging or even discharging back into the grid \citep[example case study by][]{Losel23}.
For now, we assume $0 < \chrgcurve'(\timevar)$ so that any function which is a charge curve according to \Cref{def:chargeCurve} is also a charge curve according to \Cref{def:classicCurve}.
\begin{definition}\label{def:chargeCurveMax}
    The \emph{maximum power charge curve} $\chrgcurveMax$ is the unique solution to the autonomous non-linear ordinary differential equation $\chrgcurveMax' = \ODEfunc(\chrgcurveMax)$ with boundary conditions \ref{enum:ODEleft} and \ref{enum:ODEright} of \Cref{def:chargeCurve}.
\end{definition}
$\chrgcurveMax$ is the charge curve that maximizes the charge rate throughout the recharge process.
Importantly, $\chrgcurveMax'$ gives the maximal charge rate as a function of time whereas $\ODEfunc$ gives the maximal charge rate as a function of \soc.
\Cref{fig:chargeCurveDeriv} visualizes the derivative and underlying power profile of the example charge curve from \Cref{fig:ExampleChargeCurve}.

Since we want to model the \EBSP{} with dynamic charge rates, every recharge event has to be able to choose some charge curve $\chrgcurve \in \curves(\ODEfunc)$.
Let $[\chrgStart,\chrgEnd]$ be the time window of a recharge event, $\initSoC$ the initial \soc{} at $\chrgStart$, and $\chrgcurve$ the selected charge curve.
The \soc{} at $\chrgEnd$ is then $\diffOp{\chrgcurve}(\initSoC, \chrgEnd - \chrgStart) + \initSoC$.

\begin{lemma}\label{lemma:computeChrgDiff}
    Given a finite list of time steps $(\timestep_1,\dots,\timestep_k)$, $\timestep_i > 0$, and an initial \soc{} $\initSoC\in[0,1]$, take any charge curve $\chrgcurve$ according to \Cref{def:chargeCurve} (or \Cref{def:classicCurve}) and iteratively define $\chrgvar_i = \chrgvar_{i-1} + \diffOp{\chrgcurve}(\chrgvar_{i-1}, \timestep_i)$ starting with $\chrgvar_0 = \initSoC$.
    We write $\initSoC(\chrgcurve, (\timestep_1,\dots,\timestep_k)) = \chrgvar_k$.
    Then,
    \begin{equation}
        \initSoC(\chrgcurve, (\timestep_1,\dots,\timestep_k)) = \diffOp{\chrgcurve}(\initSoC, \sum_{i=1}^k \timestep_i) + \initSoC\text{.}
    \end{equation}
\end{lemma}
\proof \cite{LBWatmos23}. \myendproof

Thus, for any equidistant time discretization of the planning horizon with step size $\timestep$, we can compute the final \soc{} of a recharge event exactly by iteratively evaluating $\chrgvar_{i} = \chrgvar_{i-1} + \diffOp{\chrgcurve}(\chrgvar_{i-1}, \timestep)$, assuming the time window aligns with the time steps.
We obtain a unary function of the \soc{} and may write $\diffOp{\chrgcurve}(\chrgvar) = \diffOp{\chrgcurve}(\chrgvar, \timestep)$.
Furthermore, let $\chrgvar(\chrgcurve, k) = \diffOp{\chrgcurve}(\chrgvar, k\timestep) + \chrgvar$ be the final \soc{} obtained after charging from $\chrgvar$ according to $\chrgcurve$ for $k$ time steps.

\begin{lemma}\label{lemma:chargeDiffDominates}
    Let $\chrgcurveMax$ be a maximum power charge curve w.r.t. a charging power profile $\ODEfunc$.
    Then for a fixed time step $\timestep > 0$, $\diffOp{\chrgcurve}(\chrgvar) \leq \diffOp{\chrgcurveMax}(\chrgvar)$ for every $\chrgvar \in [0,1]$ and $\chrgcurve \in \curves(\ODEfunc)$.
\end{lemma}
\proof Follows directly from the mean value theorem and monotonicity of $\ODEfunc$. \myendproof

\begin{corollary}
    $\chrgcurve \leq \chrgcurveMax$ for every $\chrgcurve \in \curves(\ODEfunc)$.
\end{corollary}

\begin{corollary}\label{corollary:compuationCantOverestimate}
    $\chrgvar(\chrgcurve, k) \leq \chrgvar(\chrgcurveMax, k)$ for every $\chrgvar \in [0,1]$, $\chrgcurve \in \curves(\ODEfunc)$ and $k\in\nat$.
\end{corollary}

\begin{theorem}
    Let $\initSoC$ be an initial charge state and $\timevar = k\timestep$ a charge duration with integer $k$.
    Consider the system $\chrgvar_0 = \initSoC$, $\chrgvar_{i} = \chrgvar_{i-1} + \chrgdiffvar_i$ and $0 < \chrgdiffvar_i \leq \diffOp{\chrgcurveMax}(\chrgvar_{i-1})$ for $i=1,\dots,k$.
    Then every charge curve $\chrgcurve \in \curves(\ODEfunc)$ induces a solution $(\chrgvar_0,\dots,\chrgvar_k, \chrgdiffvar_1,\dots, \chrgdiffvar_k)$ to this system and vice versa.
\end{theorem}
\proof By \Cref{lemma:computeChrgDiff}, \Cref{lemma:chargeDiffDominates} and \Cref{corollary:compuationCantOverestimate}, every charge curve $\chrgcurve \in \curves(\ODEfunc)$ yields a solution by setting $\chrgdiffvar_i = \diffOp{\chrgcurve}(\chrgvar_{i-1})$ for every $i=1,\dots,k$.
Conversely, let $(\chrgvar_0,\dots,\chrgvar_k, \chrgdiffvar_1,\dots, \chrgdiffvar_k)$ be a solution.
If $\chrgdiffvar_i = \diffOp{\chrgcurveMax}(\chrgvar_{i-1})$ for every $i=1,\dots,k$ then clearly the maximum power charge curve $\chrgcurveMax$ is induced.
Otherwise, we can construct some $\chrgcurve \in \curves(\ODEfunc)$ by appropriately throttling the charging rate $\chrgcurve'$ until $\diffOp{\chrgcurve}(\chrgvar_{i-1}) = \chrgdiffvar_i$ is achieved for every $i=1,\dots,k$.\myendproof

Note that this association is not bijective in general.
Nevertheless, given a time discretization of the planning horizon, we can introduce a \emph{charge increment variable} $\chrgdiffvar_i$ for every time step such that $\chrgvar_{i} = \chrgvar_{i-1} + \chrgdiffvar_i$ and $\chrgdiffvar_i \leq \diffOp{\chrgcurveMax}(\chrgvar_{i-1})$ to model any recharge process.
The $\chrgdiffvar_i$ imply the average applied charge rate throughout their respective time steps and thus may be incorporated into an objective function to model time-of-use prices.
Moreover, increment variables of simultaneously occurring time steps at charger slots sharing a power grid connection may be bound to adhere to the grid's maximum capacity throughout the planning horizon.
Allowing the $\chrgdiffvar_i$ to become zero can even account for a suspension of the charging process and negative values would mean a discharge back into the grid, although we will assume $\chrgdiffvar_i \geq 0$ in the following.

If we wish to use charge increment variables in mixed-integer linear programs, we have to linearize the inequality $\chrgdiffvar_i \leq \diffOp{\chrgcurveMax}(\chrgvar_{i-1})$.
\begin{definition}
    For a maximum power charge curve $\chrgcurveMax$, let $\chrgdiffDomain \defeq \setEnc{(\chrgvar, \chrgdiffvar) \in [0,1]^2 \mid \chrgdiffvar \leq \diffOp{\chrgcurveMax}(\chrgvar)}$ be the \emph{charge increment domain}.
\end{definition}

\begin{theorem}\label{theorem:convexDomain}
    Let $\chrgcurveMax$ be a maximum power charge curve w.r.t. a charging power profile $\ODEfunc$. Furthermore, let $\ODEfunc$ be concave. Then the corresponding charge increment domain $\chrgdiffDomain$ is convex for any time step size $\timestep > 0$.
\end{theorem}
\proof \cite{LBWatmos23}. \myendproof

\begin{figure}[ht]
    {\includegraphics[width=\textwidth]{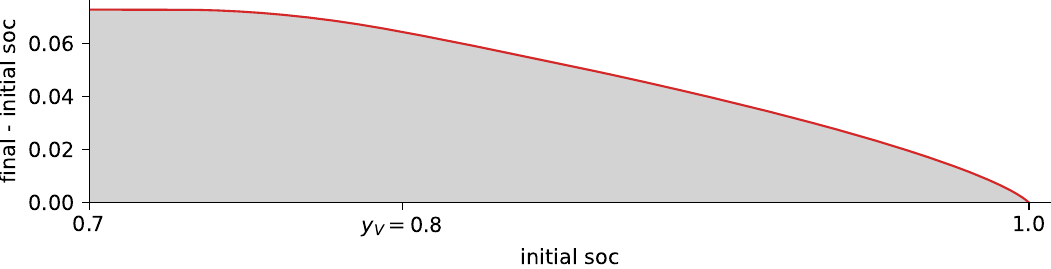}}
    \caption{The shaded area is the charge increment domain $\chrgdiffDomain$ on $[0.7, 1]$ for the charge curve presented in \Cref{fig:ExampleChargeCurve} with $\timestep = 5\text{min}$.}
    \label{fig:ExampleChargeDomain}
\end{figure}

Observe that the charging power profile of the example curve in \Cref{fig:chargeCurveDeriv} is indeed concave.
A corresponding convex charge increment domain is shown in \Cref{fig:ExampleChargeDomain}.
We therefore suggest, and this is our main modeling contribution, to interpolate the function graph of $\diffOp{\chrgcurveMax}$ (as a unary function of the \soc) with a linear spline $\linSpline{\diffOp{\chrgcurveMax}}$ and to replace $\chrgdiffvar_i \leq \diffOp{\chrgcurveMax}(\chrgvar_{i-1})$ by inequalities of the form $\chrgdiffvar_i \leq \slope_j\chrgvar_{i-1} + \offset_j$ for each linear segment $j$.
Interpolating the charge increment function instead of the charge curve has many advantages.
If $\ODEfunc$ (or alternatively $\chrgcurveMax'$ parametrized as a function of the \soc) is concave, then this interpolation is guaranteed to underestimate the exact boundary by \Cref{theorem:convexDomain}.
Importantly, $\chrgdiffvar_i \leq \diffOp{\chrgcurveMax}(\chrgvar_{i-1})$ is a valid inequality for the set of feasible charge increment variable values and we simply restrict the set of available charge curves to $\approxOp{\curves}(\ODEfunc) \subset \curves(\ODEfunc)$.
Any sequence of increment variables that lie on the interpolation spline induce a maximum power charge curve $\approxOp{\chrgcurveMax}$ on $\approxOp{\curves}(\ODEfunc)$.
This $\approxOp{\chrgcurveMax}$ is then an approximation of $\chrgcurveMax$ which satisfies $\diffOp{\approxOp{\chrgcurveMax}} \leq \diffOp{\chrgcurveMax}$ by \Cref{lemma:chargeDiffDominates}, so the bounds from \Cref{sec:rechargeModeling} hold.

\begin{proposition}
    \label{prop:incDomainError}
    Let $\linSpline{\diffOp{\chrgcurveMax}}$ be a piecewise linear spline interpolation of $\diffOp{\chrgcurveMax}$. Then the approximation error is
    \begin{equation}
        \norm{\linSpline{\diffOp{\chrgcurveMax}} - \diffOp{\chrgcurveMax}} \leq \frac{\timestep h^2}{8} \norm{\ODEfuncCV''}
    \end{equation}
    where $h$ is the width of the largest linear segment.
\end{proposition}
\proof \cite{LBWatmos23}. \myendproof

Selecting $\chrgvar_1 = \chrgcurveMax(\CCCVbreakTime - \timestep)$ as the first inner interpolation point, we know that the linear spline is exact on the first segment and can deduce that $h \leq 1 - \CCCVbreakSoC + \ODEfuncCC\timestep$, which is the width of the increment domain segment corresponding to the CV-phase.

\begin{corollary}\label{cor:CIDerrorBound}
    Using a linearized charge increment domain satisfying \Cref{theorem:convexDomain}, we have
    \begin{equation}
        \abs{\errorFunc(\course_j)} \leq \numChrgEvents(j) \frac{\timestep \enc{1 - \CCCVbreakSoC + \ODEfuncCC\timestep}^2}{8} \norm{\ODEfuncCV''} \in o(\timestep^3)
    \end{equation}
    along a vehicle course $\course$. Recall that $\numChrgEvents(j)$ gives the number of recharge events occurring before $C_j$.
\end{corollary}

Note that the error can be controlled both by the interpolation step-size as well as the step size of the time discretization, and decreases by order 2 in the first and even by order 3 in the second case. 
Moreover, the error estimate provides a bound for the maximal underestimation of the exact charge state.
It will turn out in the empirical evaluation in \Cref{sec:approxErrorSpan} that this  approximation scheme is quasi exact on real-life instances if reasonable discretization step sizes are used.

Should $\chrgdiffDomain$ not be convex, we can still attempt to fit a concave piecewise linear function to its boundary, but then we either lose $\chrgdiffvar_i \leq \diffOp{\chrgcurveMax}(\chrgvar_{i-1})$ being a valid inequality or we have to be content with a larger approximation error.

\section{Electric Bus Scheduling with the Charge Increment Domain}\label{sec:MIP}
We explain our mixed-integer linear program for the \EBSP{} here.
It combines the new charging model with a standard resource constrained binary multi-commodity flow formulation on a directed graph $\DAG=(\nodes, \arcs)$ with a node for every depot in $\deps$ and trip in $\trips$.
As in other models, the entire service horizon is discretized using equidistant time steps and we introduce a \emph{timeline} for every charger slot, i.e., we create a node $\chrgslot_i$ for every slot $\chrgslot\in\chrgslots$ and time event $i = 0,\dots,\horizonEnd$.
This is so we can incorporate the charge increment domain linearization and track charge slot occupation.
There may be additional auxiliary nodes and timelines to, for example, allow intermittent parking at the depot.

Every two consecutive charge nodes $\chrgslot_{i-1}$ and $\chrgslot_i$ on the timeline of a charge slot $\chrgslot$ are connected by a \emph{recharge arc} $\arc(\chrgslot, i) \defeq (\chrgslot_{i-1}, \chrgslot_i)$.
We denote the set of all charge nodes by $\chrgNodes$, the set of all recharge arcs by $\rechargeArcs$, and the recharge arcs belonging to a particular charge slot by $\rechargeArcs(\chrgslot)$.
The remaining arcs in $\arcs$ encode feasible connections between the trip nodes, depot nodes, and charge slot timeline nodes.
$\depotPullouts$ denotes the set of depot pull-out arcs, i.e., arcs with a depot node as their tail.

We have a \emph{plan type} $\commodity$ for every combination of vehicle type $\vcl\in\vcls$ and depot $\dep\in\deps$ that may house it.
We denote the set of all plan types by $\commodities$ and the subset belonging to electric vehicle types by $\commoditiesE$.
Not every plan type may be admissible on every arc, so $\commodities(\arc)$ denotes the set of plan types allowed on arc $\arc$.
We can use this to enforce that buses terminate their service at their depot of origin by restricting plan types on depot adjacent arcs and block non-electric types from accessing charge stations.

As such, our MILP has binary flow variables $\flowvar_{\arc}^{\commodity}$ per plan type on each admissible arc, which are subject to usual flow conservation constraints per commodity at every node
\begin{equation}
    \sum_{\substack{\arc\in\inset(\node): \\ \commodities(\arc)\ni \commodity}} \flowvar_{\arc}^{\commodity} = \sum_{\substack{\arc\in\outset(\node): \\ \commodities(\arc)\ni \commodity}} \flowvar_{\arc}^{\commodity} \quad\forall\, \node \in \nodes \setminus \deps, ~\commodity\in\commodities\text{,}\label{MIP:flowConservation}
\end{equation}
where $\inset(\node)$ denotes the set of incoming and $\outset(\node)$ the set of outgoing arcs of node $\node$.
We require that every trip is serviced by exactly one bus
\begin{equation}
    \sum_{\substack{\arc\in\outset(\trip), \\ \commodity\in\commodities(\arc)}}\flowvar_{\arc}^{\commodity} = 1 \quad\forall \trip \in \trips \label{MIP:tripDemand}
\end{equation}
and every other non-depot node may be used by at most one bus
\begin{equation}
    \sum_{\substack{\arc\in\outset(\node), \\ \commodity\in\commodities(\arc)}} \flowvar_{\arc}^{\commodity} \leq 1 \quad\forall \node\in\nodes\setminus(\deps \cup \trips)\text{.} \label{MIP:outCapacity}
\end{equation}
Moreover, there may be so-called \emph{vehicle-mix constraints} bounding admissible fleet compositions and their assignment to depots.
This, for example, accounts for a depot's total overnight parking space in relation to different vehicle type sizes.
They are given by $\mixCons \subset 2^{\commodities}$, a lower and upper bound $\mixLower_{\mixTypes}$ and $\mixUpper_{\mixTypes}$ per $\mixTypes \in \mixCons$, and coefficients $\mixCoeff_{\mixTypes}^{\commodity}$ for $\commodity \in \mixTypes \in \mixCons$.
The mix-constraints bound vehicle course assignments by
\begin{equation}
    \mixLower_{\mixTypes} \leq \sum_{(\vcl, \dep) \in \mixTypes} \mixCoeff_{\mixTypes}^{(\vcl, \dep)} \sum_{\arc\in\outset(\dep)} \flowvar_{\arc}^{(\vcl, \dep)} \leq \mixUpper_{\mixTypes} \quad \forall \mixTypes \in \mixCons \text{.} \label{MIP:mixConstraint}
\end{equation}

Variables $\chrgvar_{\arc}$ give the remaining driving range of the active bus on the arc just after the tail node.
Tracking \soc{} on the arcs has been proposed by \cite{FMJL19}.
We normalize the battery capacities and energy requirements to $[0,1]$ to account for different battery capacities of different bus types.
The energy variables are then coupled to the arc's total flow by
\begin{equation}
    \sum_{\commodity\in\commoditiesE(\arc)} \flowvar_{\arc}^{\commodity} \geq \chrgvar_{\arc} \quad \forall \arc \in \arcs \setminus \depotPullouts \label{MIP:chrgUpperBound}
\end{equation}
and we enforce that buses leave the depot fully charged by stipulating
\begin{equation}
    \sum_{\commodity\in\commoditiesE(\arc)} \flowvar_{\arc}^{\commodity} = \chrgvar_{\arc} \quad \forall \arc \in \depotPullouts\text{.} \label{MIP:chrgPullout}
\end{equation}
The remaining driving range is propagated as an energy flow via
\begin{equation}
    \sum_{\substack{\arc\in\inset(\node), \\ \commodity \in \commoditiesE(\arc)}}\consumption_{\arc}^{\commodity} \flowvar_{\arc}^{\commodity} = \sum_{\arc\in\inset(\node)} \chrgvar_{\arc} - \sum_{\arc\in\outset(\node)} \chrgvar_{\arc} \quad \forall \node \in \bar{\nodes} \label{MIP:eneryFlow}
\end{equation}
where $\bar{\nodes} \defeq \nodes \setminus (\deps \cup \chrgNodes)$ and the $\consumption_{\arc}^{\commodity}$ give the relative range requirement for a bus of plan type $\commodity$ to drive along arc $\arc$ and its head node.
The remaining range just after node $\node$ is the remaining range at the beginning of the active incoming arc minus the corresponding consumption.
Non-electric bus types consequently incur an energy consumption of zero.

For charge nodes, we have to include the charge increment variables on the recharge arc, so the energy propagation constraint becomes
\begin{equation}
    \sum_{\substack{\arc\in\inset(\chrgslot_i), \\ \commodity \in \commoditiesE(\arc)}}\consumption_{\arc}^{\commodity} \flowvar_{\arc}^{\commodity} - \sum_{\commodity \in \commoditiesE(\arc(\chrgslot, i))} \chrgdiffvar_{\arc(\chrgslot, i)}^{\commodity} \\
    = \sum_{\arc\in\inset(\chrgslot_i)} \chrgvar_{\arc} - \sum_{\arc\in\outset(\chrgslot_i)} \chrgvar_{\arc} \quad \forall \arc(\chrgslot, i) \in \rechargeArcs\text{,} \label{MIP:energyFlowCharging}
\end{equation}
with an energy consumption of usually zero along the recharge arc.
In our model, the charge increment domain is described by a piecewise linear function. Let $\numLinSegments(\chrgslot, \commodity)$ be the number of linear segments for charger $\chrgslot$ and admissible plan type $\commodity$, and let $\slope_j(\chrgslot, \commodity)$ denote the $j$th slope and $\offset_j(\chrgslot, \commodity)$ the corresponding offset.
The underlying increment domain shall be convex, so $0 = \slope_1(\chrgslot, \commodity) > \slope_2(\chrgslot, \commodity) > \dots > \slope_{\numLinSegments(\chrgslot, \commodity)}(\chrgslot, \commodity)$ must hold.
We can then easily handle recharge processes with constraints
\begin{equation}
    \offset_1(\chrgslot, \commodity) \flowvar_{\arc(\chrgslot, i)}^{\commodity} \geq \chrgdiffvar_{\arc(\chrgslot, i)}^{\commodity} \quad \substack{\forall \arc(\chrgslot, i) \in \rechargeArcs, \\ \commodity \in \commoditiesE(\arc(\chrgslot, i))} \label{MIP:chrgIncCoupling}
\end{equation}
and
\begin{equation}
    \slope_j(\chrgslot, \commodity)\, \chrgvar_{\arc(\chrgslot, i)} + \offset_j(\chrgslot, \commodity) \geq \chrgdiffvar_{\arc(\chrgslot, i)}^{\commodity} \quad \substack{\forall \arc(\chrgslot, i) \in \rechargeArcs, \\ \commodity \in \commoditiesE(\arc(\chrgslot, i)), \\ j = 2, \dots, \numLinSegments(\chrgslot, \commodity)} \label{MIP:chrgIncDomain}
\end{equation}
to bound the increment variable according to both its recharge arc being active and the initial energy state at the beginning of the recharge time step.
The maximum grid load is also easily enforced as
\begin{equation}
    \sum_{\substack{\chrgslot \in \accessPoint, \\ \commodity \in \commoditiesE(\arc(\chrgslot, i))}} \maxPowerCoeff(\chrgslot, \commodity)\, \chrgdiffvar_{\arc(\chrgslot, i)}^{\commodity} \leq \maxPowerDraw_i(\accessPoint) \quad \forall \accessPoint \in \gridAccess,~ i \in [\horizonEnd]\text{,} \label{MIP:gridCapacity}
\end{equation}
where $\gridAccess$ is a partition of the charge slots $\chrgslots$ by shared power grid access.
For each such access point $\accessPoint \in \gridAccess$, there is a maximum amount of available power $\maxPowerDraw_i(\accessPoint)$ per time step.
Due to the normalization of the battery capacities, we need coefficients $\maxPowerCoeff(\chrgslot, \commodity)$ to scale the charge increment variable values accordingly.

Every plan type admits operational costs $\cost_{\arc}^{\commodity}$ per arc.
Additionally, there may be time-of-day dependent electricity costs $\costtou_{\arc_i}^{\commodity}$ for the charge increment variables on recharge arcs.
The full MILP is 
\begin{align}
    && \min && \sum_{\substack{\arc\in\arcs, \\ \commodity\in\commodities(\arc)}} \cost_{\arc}^{\commodity} \flowvar_{\arc}^{\commodity} & + \sum_{\substack{\arc \in \rechargeArcs, \\ \commodity\in\commoditiesE(\arc)}} \costtou_{\arc}^{\commodity} \chrgdiffvar_{\arc}^{\commodity} && \label{MIP:objective} \\
    && \text{s.t.} && \text{\eqref{MIP:flowConservation} } & \text{to \eqref{MIP:gridCapacity}} && \\
    &&&& \flowvar_{\arc}^{\commodity} & \in \setEnc{0,1} && \forall \arc\in\arcs,\commodity\in\commodities(\arc) \label{PBSP:outCapacity} \\
    &&&& \chrgvar_{\arc} & \geq 0 && \forall \arc \in \arcs \label{MIP:chrgvarDomain} \\
    &&&& \chrgdiffvar_{\arc}^{\commodity} & \geq 0 && \forall \arc \in \rechargeArcs, \commodity \in \commoditiesE(\arc) \label{MIP:chrgIncVarDomain}
\end{align}

Note that if we replace \eqref{MIP:chrgIncCoupling} and \eqref{MIP:chrgIncDomain} by
\begin{equation}
    \diffOp{\chrgcurveMax(\chrgslot, \commodity)}(\chrgvar_{\arc(\chrgslot, i)}) \,\flowvar_{\arc(\chrgslot, i)}^{\commodity} \geq \chrgdiffvar_{\arc(\chrgslot, i)}^{\commodity} \quad \substack{\forall \arc(\chrgslot, i) \in \rechargeArcs, \\ \commodity \in \commoditiesE(\arc(\chrgslot, i))}
\end{equation}
where $\chrgcurveMax(\chrgslot, \commodity)$ is the maximum power charge curve of the corresponding charger and vehicle type combination, the model becomes exact but in general non-linear.
In both cases, our approach extends a standard model by a charging component of controlled accuracy, and the same extension works for any similar model. 

Some battery types have to be preconditioned or require other preparatory work before they can be recharged despite occupying a charge slot.
To account for this, additional constraints of the form $\offset_1(\chrgslot, \commodity) \flowvar_{\arc(\chrgslot, j)}^{\commodity} \geq \chrgdiffvar_{\arc(\chrgslot, i)}^{\commodity}$ with appropriate $j < i$ can be added.

\subsection{Strengthening Inequalities}\label{sec:ineq}
It is straightforward to tighten inequalities \eqref{MIP:chrgUpperBound} and $\chrgvar \geq 0$ for the LP relaxation.
It is clear that the maximum and minimum energy state at the head of every arc is dictated by the surrounding network, i.e., there is a minimum amount of energy required to reach and to leave a node.
As such, for every node $\node$ let $\inPaths(\node, \commodity)$ denote the set of all paths from a depot or charge node to $\node$ using arcs admissible for plan type $\commodity$, and let $\outPaths(\node, \commodity)$ denote the corresponding set of paths from $\node$ to a depot or charge node.
Defining
\begin{equation}
    \maxChrg_{\node}^{\commodity} \defeq \max\setEnc{1 - \sum_{\arc \in \arcPath}\consumption_{\arc}^{\commodity} \smid \arcPath \in \inPaths(\node, \commodity)}
\end{equation}
and
\begin{equation}
    \minReqChrg_{\node}^{\commodity} \defeq \min\setEnc{\sum_{\arc \in \arcPath}\consumption_{\arc}^{\commodity} \smid \arcPath \in \outPaths(\node, \commodity)}\text{,}
\end{equation}
we can replace \eqref{MIP:chrgUpperBound} and $\chrgvar \geq 0$ by
\begin{equation}
    \sum_{\commodity\in\commoditiesE((\node_t, \node_h))} (\consumption_{(\node_t, \node_h)}^{\commodity} + \minReqChrg_{\node_h}^{\commodity}) \flowvar_{(\node_t, \node_h)}^{\commodity} \leq \chrgvar_{(\node_t, \node_h)} \quad \forall (\node_t, \node_h) \in \arcs \setminus \depotPullouts
\end{equation}
and
\begin{equation}
    \chrgvar_{(\node_t, \node_h)} \leq \sum_{\commodity\in\commoditiesE((\node_t, \node_h))} \maxChrg_{\node_t}^{\commodity} \flowvar_{(\node_t, \node_h)}^{\commodity} \quad \forall (\node_t, \node_h) \in \arcs \setminus \depotPullouts
\end{equation}
in the MILP.
We have found that these inequalities do separate optimal LP solutions to the formulation while requiring only trivial preprocessing computations.

\section{Computational Study}\label{sec:study}
We tested our model on sixteen anonymized real-life \EBSP{} instances (see \Cref{tab:instances}) using Gurobi 11.0.0 \citep{gurobi} on a heterogeneous computation cluster, where each experiment was bound to four physical cores and eight threads to simulate execution on an ordinary desktop workstation.

\begin{table}[htb]
    {\begin{tabular*}{\textwidth}{@{\extracolsep{\fill}} c | *{7}{r}}
        \toprule
        instance & \makecell[c]{electric\\bus types} & \makecell[c]{non-electric\\bus types} & \makecell[c]{depots} & \makecell[c]{charge\\slots} & \makecell[c]{grid access\\points} & \makecell[c]{trips} & \makecell[c]{deadheads} \\
        \midrule
        A & 1 & 0 & 1 & 3 & 1 & 121 & 992 \\
        B & 1 & 0 & 1 & 2 & 1 & 123 & 1\,078 \\
        C & 2 & 0 & 1 & 3 & 1 & 146 & 3\,126 \\
        D & 1 & 0 & 1 & 3 & 1 & 185 & 1\,769 \\
        E & 1 & 0 & 1 & 8 & 1 & 189 & 1\,977 \\
        F & 1 & 0 & 1 & 3 & 1 & 232 & 1\,484 \\
        G & 1 & 0 & 1 & 5 & 2 & 232 & 2\,064 \\
        H & 1 & 1 & 1 & 6 & 1 & 333 & 6\,363 \\
        I & 1 & 1 & 1 & 7 & 2 & 333 & 7\,859 \\
        J & 1 & 0 & 1 & 14 & 1 & 678 & 15\,589 \\
        K & 1 & 0 & 1 & 43 & 1 & 709 & 14\,431 \\
        L & 1 & 0 & 1 & 37 & 1 & 709 & 17\,779 \\
        M & 1 & 0 & 1 & 34 & 1 & 709 & 21\,343 \\
        N & 2 & 1 & 2 & 12 & 1 & 822 & 12\,390 \\
        O & 1 & 1 & 1 & 10 & 1 & 837 & 111\,590 \\
        P & 1 & 0 & 1 & 28 & 1 & 1\,207 & 42\,610 \\
        \bottomrule
    \end{tabular*}}
    \caption{The number of electric bus types, non-electric bus types, depots, charge slots, grid access points, timetables passenger trips, and deadhead connections for every instance.
    The number of charge slots give the maximum number of buses that can charge simultaneously.
    They are distributed over one or two power grid access points.
    The instances are in ascending order with respect to the number of trips.}
    \label{tab:instances}
\end{table}

Instance G extends F, and I extends H by an additional opportunity fast-charging terminal.
Instances K, L and M are variations of the same network with different charging technology and bus types.
For the number of deadheads note that generally only short and direct connections are explicitly given.
Long deadheads including a stopover at a parking facility or depot are given implicitly by parking timelines and connecting pull-in and pull-out deadheads.

\subsection{Discretization Parameter Analysis}\label{sec:discParam}
Our model has two discretization parameters: The time step size $\timestep$ along charge slot timelines and the number of linear segments $\numLinSegments$ for the charge increment domain.
We want to determine the best parameters empirically and sampled $\numLinSegments=2,3,4,10$ and $\timestep$ as $60$, $300$ and $600$ seconds as a priori reasonable choices by solving the MILP formulation from \Cref{sec:MIP} without grid capacity limits \eqref{MIP:gridCapacity} on our test instances.

To demonstrate that exact charge models are needed, we also pre-computed vehicle schedules under a fully linear approximation of the respective charge curve.
We used the method described in \citep{BLLW24} followed by a further Gurobi optimization using our MILP formulation, where we dropped constraints \eqref{MIP:chrgIncDomain} and selected a time step size of $60$ seconds.
This solved instances A, C, G, J, L, O, and P to (almost) optimality within numerical tolerances under the linear charging model.
The resulting schedules were then evaluated as initial solution candidates and we analyzed their energy-feasibility under our more exact model using constraints \eqref{MIP:chrgIncDomain}.
We further cold-started every configuration and recorded the best obtained solution and lower bound after one and twelve hours of runtime.

The numerical results presented in \Crefrange{tab:discParam1}{tab:discParam8} in the appendix show that for ten of the sixteen instances the initial vehicle schedule is not energy-feasible under the more exact charging model for all tested parameter.
On five additional instances the initial vehicle schedule fails for some of the parameter configurations, although on G, J, L, and P this is likely due to the coarser time discretization.
On O and P the initial schedule is infeasible for $\numLinSegments=2$ but feasible for $\numLinSegments > 2$, which indicates interference by an approximation error.
In both cases, Gurobi manages to repair the solution after pre-processing.

We observe that the lower bounds on A are slightly larger than the objective value of the optimal schedule under linear charging for every configuration, i.e., we have a non-zero gap between the charging models.
On C, G, J, L, O, and P we can obtain energy-feasible schedules whose objective is equal to the optimal reference solution within numerical tolerances.
On the remaining instances other than K we can find solutions whose objective value has the same order of magnitude as the reference solution.
Moreover, the lower bound for M is larger than the objective value of the reference solution for $\timestep = 600$ seconds only, indicating that a time step size of ten minutes is too coarse.

\begin{figure}[htb]
    \includegraphics[width=\textwidth]{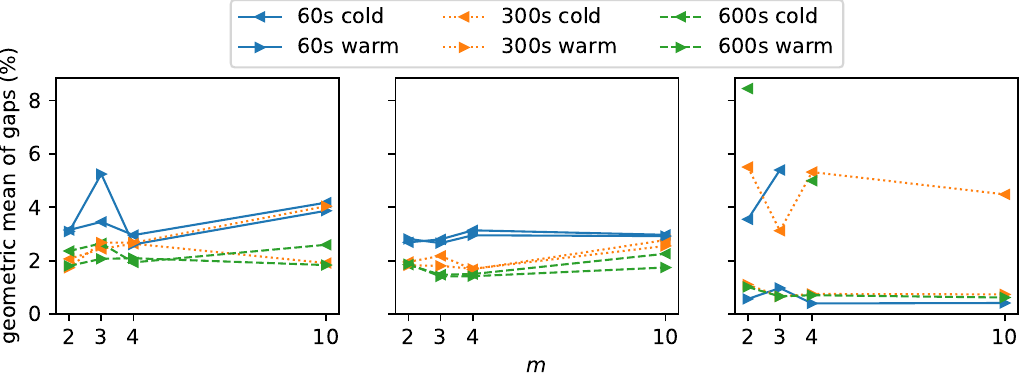}
    \caption{The geometric means of the gaps for instances $Q_1 = \setEnc{A, B, D, E, F, G, J, N}$ on the left, $Q_2 = \setEnc{C, K, M, O}$ in the middle, and $Q_3 = \setEnc{H, I, L, P}$ on the right.
    The means are plotted over the number of linear segments of the charge increment domain $\numLinSegments$ per time step size $\timestep$ and whether an initial solution was provided (``warm") or not (``cold").
    Missing markers in the plot on the right indicate that this particular configuration failed to produce a feasible solution on all four instances in the batch within twelve hours.}
    \label{fig:discParamGaps}
\end{figure}

Finally, to help determine the best parameters, \Cref{fig:discParamGaps} presents the geometric mean over the gaps of each tested configuration.
We partition our sixteen instances into three batches according to their apparent difficulty for Gurobi: $Q_1 = \setEnc{A, B, D, E, F, G, J, N}$, the instances which have a feasible solution for all configurations after one hour of runtime, $Q_2 = \setEnc{C, K, M, O}$, the instances not already in $Q_1$ which have a feasible solution for all configurations after twelve hours of runtime, and $Q_3 = \setEnc{H, I, L, P}$, the remaining instances.
If a configuration produced no feasible solution, the instance did not contribute a gap to the mean for that parameter choice for the plot on the right.

From the analysis of the numerical results we conclude that we should choose $\numLinSegments > 2$ and $\timestep < 600$ seconds.
Then, $\timestep = 300$ seconds clearly outperforms $\timestep = 60$ seconds in \Cref{fig:discParamGaps}.
Note that \citeauthor{VLD23}, who also use discrete timelines at charge slots, settled on a time step size of five minutes after their sensitivity analysis as well.
For the number of linear segments $\numLinSegments$, $3$ and $4$ have a roughly comparable performance.
Note that the $300$ seconds ``cold'' plot on the right of \Cref{fig:discParamGaps} is an artifact caused by the low availability of solutions in that batch.
Since the approximation error shrinks with the square of the width of the interpolation segments by \Cref{prop:incDomainError}, we select $\numLinSegments = 4$ for further experiments.

\subsection{Over- and Underestimating the Increment Domain}\label{sec:approxErrorSpan}
\begin{figure}[htb]
    \includegraphics[width=\textwidth]{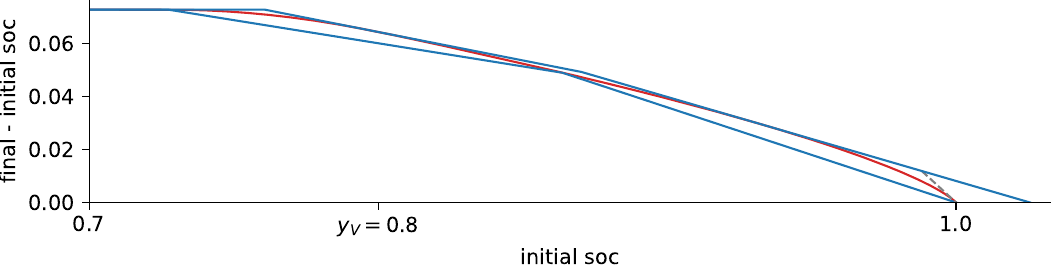}
    \caption{The increment domain section from \Cref{fig:ExampleChargeDomain} with an example underestimating linear spline interpolation and an overestimating linear spline of tangential lines using $\numLinSegments = 3$ segments.
    $\numLinSegments = 4$ is already too tight to properly see the difference between both splines and the exact domain.
    Note that there is an implicit bound of $\chrgdiffvar \leq 1 - \chrgvar$ in our model which is represented by the dashed line segment near $1.0$ on the initial \soc{} axis.}
    \label{fig:ChargeDomainApprox}
\end{figure}

To gauge the impact of the approximation error of our charge behavior linearization in practice, we solved our instance set with an overestimation of the increment domain.
Instead of constructing constraints \eqref{MIP:chrgIncCoupling} and $\eqref{MIP:chrgIncDomain}$ from a linear interpolation on equidistant grid points $(\chrgvar_i, \diffOp{\chrgcurveMax}(\chrgvar_i))$, we used the tangential lines of $\diffOp{\chrgcurveMax}$ at $(\chrgvar_i + \chrgvar_{i+1}) / 2$, yielding a piecewise linear, overestimating approximation, see \Cref{fig:ChargeDomainApprox}.

\begin{figure}[htb]
    \includegraphics[width=\textwidth]{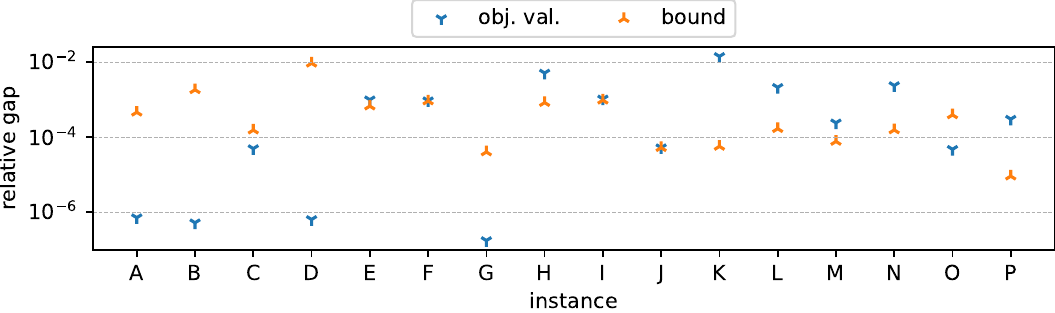}
    \caption{Relative gaps between final objective values and bounds of under- and overestimating the increment domain boundary.
    The instances were solved for twelve hours using $\numLinSegments=4$ and $\timestep = 300s$.
    The runs were warm started as in the parameter analysis experiment.
    The gaps were computed by $\abs{\overline{b} - \underline{b}} / ((\overline{b} + \underline{b})/2)$ where $\overline{b}$ is the objective value or lower bound of the overestimation and $\underline{b}$ of the underestimation.}
    \label{fig:overUnderEst}
\end{figure}

As can be deduced from \Cref{fig:overUnderEst}, there is no significant approximation error on our instance test set in the sense that both the under- and overestimation produce numerically equivalent objective values and bounds.
Note that we have three pools of charge curves satisfying $\approxOp{\curves}(\ODEfunc) \subset \curves(\ODEfunc) \subset \approxOp{\curves}'(\ODEfunc)$ in this setting such that the exact maximum power charge curve sits in between the maximum power curve of the under- and of the overestimation.
We can therefore be confident that our model is sufficiently accurate on real instances.

\subsection{Limiting the Grid Load}\label{sec:gridLoadTest}
\begin{figure}[htb]
    \includegraphics[width=\textwidth]{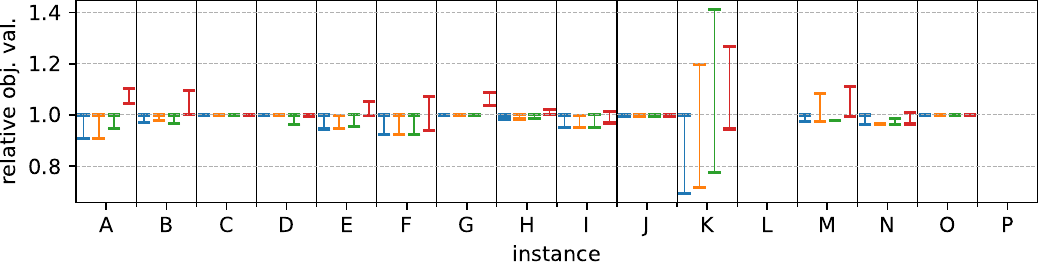}
    \caption{Final relative objective values and lower bounds after solving the instances for twenty four hours using $\numLinSegments=4$ and $\timestep = 300s$ and limiting the peak power grid load to $1.0$, $0.75$, $0.5$ and $0.25$ times the peak load from the solution obtained in \Cref{sec:discParam}.
    The values are printed from left to right per instance under decreasing peak load limit.
    The marker for the objective and lower bound are connected by a vertical line as a visual aid.
    Gurobi failed to find any incumbent for instances L and P, just like for the cold-started runs in \Cref{sec:discParam}.
    Also maintaining previously seen behavior, instance K is left with large integrality gaps.}
    \label{fig:gridLimit}
\end{figure}

For the last experiment, we are testing the impact of limiting the grid load by imposing all constraints of type \eqref{MIP:gridCapacity}, which were absent for the previous experiments.
We set the admissible peak load to $1.0$, $0.75$, $0.5$ and $0.25$ times the peak power requirement of the respective best solution for $\numLinSegments = 4$ and $\timestep = 300s$ from \Cref{sec:discParam}.
Since a given vehicle schedule may become infeasible under tighter peak load limits, the computations were cold-started to keep the runs comparable between the different load limits.
To compensate for the expected increased computational complexity, the experiment ran for 24 hours.
The results are presented in \Cref{fig:gridLimit} where we normalized all objective values and lower bounds such that the final objective value obtained under the full maximum peak load is at $1.0$ for every instance.

The results clearly indicate that up until half the maximum peak, the charging load can be distributed more evenly over the planning horizon without a significant cost penalty, but a quarter maximum peak may cut off the unrestricted reference solution.
We take a closer look at the grid load distribution for two of the instances to demonstrate so-called peak shaving.

\begin{figure}[!htb]
    \includegraphics[width=\textwidth]{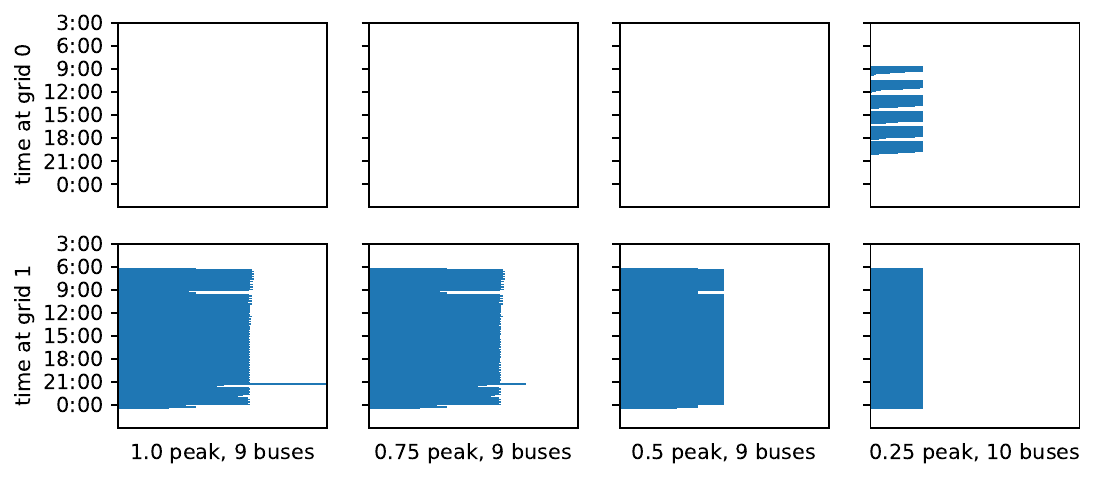}
    \caption{Total grid load per time step relative to the total peak load of instance G.
    There are two access points, one at the depot (grid 0) and one for a fast charger at a central location of turn points (grid 1).}
    \label{fig:powerReqG}
\end{figure}

Instance G (\Cref{fig:powerReqG}) possesses slow chargers at the depot and fast opportunity chargers at a separate location.
In the original solution, there is a constant base load of a single bus getting topped-off. 
Regularly, there is some overlap between the charge events and two buses are drawing power, although not necessarily at the highest possible rate.
This causes the fuzzy boundary in the middle and on the right-hand-side of the plots for $1.0$ and $0.75$.
At half past nine in the evening, a third bus is charged simultaneously, causing a load peak.
For the $0.5$ limit, we can clearly observe a more even distribution of the total grid load.
Going even lower requires the depot charger to be used as well and the added downtime forces another bus to be deployed.
The lower bound in \Cref{fig:gridLimit} indicates that this is indeed the optimal fleet size under the tighter load limit.
Of course, on this particular example we could have dealt with the undesired load peak in the evening by reducing the number of charger slots to two, but peak shaving is not as easy on larger instances. 

\begin{figure}[!htb]
    \includegraphics[width=\textwidth]{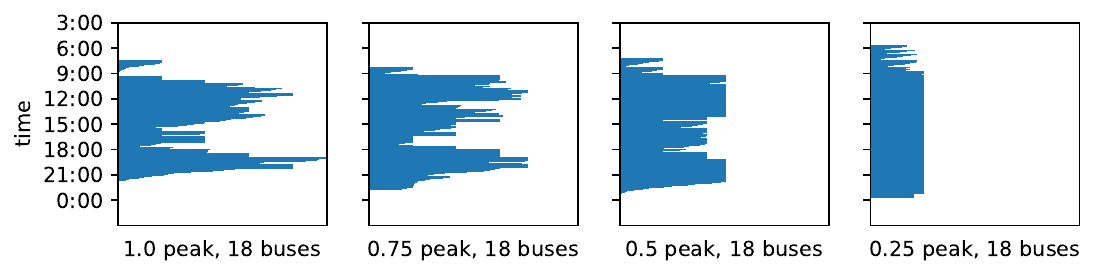}
    \caption{Total realtive grid load per time step relative to the total peak load of instance H.}
    \label{fig:powerReqH}
\end{figure}

Public transport timetables are often tightened during the morning and afternoon to accommodate school traffic, which leads to charging event peaks at midday and in the evening.
Instance H exhibits this typical behavior in \Cref{fig:powerReqH}.
Consulting \Cref{fig:gridLimit} as well we see that we can shave these grid load peaks at only minor cost penalties.
We could not have replicated this result by limiting the number of available charger slots because all four solutions have the instance's six chargers fully occupied at some point during the evening.
The buses are instead recharged at a throttled rate, just as we intended with the model.

\section{Conclusion}
We have given a brief overview of how the literature deals with the non-linear charging behavior for the electric bus scheduling problem and examined the impact of approximation errors on solutions.
We have motivated and developed a novel model formulation for \EBSP{} with a well understood theoretical approximation behavior, that is based on a piecewise linear interpolation of the charge increment curve. It extends standard bus scheduling models in a simple but flexible way, and allows to take current and emerging operator requirements into account, in particular dynamic charging and power grid load concerns.
Finally, we verified the usefulness of this model on a test set of real-life instances.

~

\acknowledgement{This work has been conducted within the Research Campus MODAL funded by the German Federal Ministry of Education and Research (BMBF) (fund number 05M20ZBM).}

\bibliography{bib.bib}

\begin{thebibliography}{}

\bibitem[Abdelwahed et~al., 2020]{ABBCK20}
Abdelwahed, A., van~den Berg, P., Brandt, T., Collins, J., and Ketter, W.
  (2020).
\newblock {Evaluating and Optimizing Opportunity Fast-Charging Schedules in
  Transit Battery Electric Bus Networks}.
\newblock {\em Transportation Science}, pages 1--15.

\bibitem[Auberle and Reidelst{\"u}rz, 2024]{AR24}
Auberle, A. and Reidelst{\"u}rz, M. (2024).
\newblock {Elektrifizierung der Busflotte: V{\"o}lklinger Verkehrsbetriebe
  verlassen sich auf Simulationen}.
\newblock {\em Der Nahverkehr Elektrobus-Spezial 2024}, 42:26–28.
\newblock (German).

\bibitem[Bagherinezhad et~al., 2020]{BPLP20}
Bagherinezhad, A., Palomino, A.~D., Li, B., and Parvania, M. (2020).
\newblock {Spatio-Temporal Electric Bus Charging Optimization With Transit
  Network Constraints}.
\newblock {\em IEEE Transactions on Industry Applications}, 56(5):5741--5749.

\bibitem[Bertossi et~al., 1987]{BCG87}
Bertossi, A.~A., Carraresi, P., and Gallo, G. (1987).
\newblock {On Some Matching Problems Arising in Vehicle Scheduling Models}.
\newblock {\em Networks}, 17:271--281.

\bibitem[Bornd\"{o}rfer et~al., 2024]{BLLW24}
Bornd\"{o}rfer, R., L\"{o}bel, A., L\"{o}bel, F., and Weider, S. (2024).
\newblock {Solving the Electric Bus Scheduling Problem by an Integrated Flow
  and Set Partitioning Approach}.
\newblock Technical report, Zuse Institute Berlin.
\newblock not yet published, submitted for review.

\bibitem[Bunte and Kliewer, 2010]{BK10}
Bunte, S. and Kliewer, N. (2010).
\newblock {An Overview on Vehicle Scheduling Models}.
\newblock {\em Public Transport}, 1:299--317.

\bibitem[de~Vos et~al., 2023]{VLD23}
de~Vos, M.~H., van Lieshout, R.~N., and Dollevoet, T. (2023).
\newblock {Electric Vehicle Scheduling in Public Transit with Capacitated
  Charging Stations}.
\newblock {\em Transportation Science}, 2023:1--16.

\bibitem[Diefenbach et~al., 2023]{DEG23}
Diefenbach, H., Emde, S., and Glock, C.~H. (2023).
\newblock Multi-depot electric vehicle scheduling in in-plant production
  logistics considering non-linear charging models.
\newblock {\em European Journal of Operational Research}, 306(2):828--848.

\bibitem[Dietmannsberger, 2023]{Hamburg23}
Dietmannsberger, M. (2023).
\newblock {100 Busse geschafft - und jetzt? Hamburger Wege zur Verbesserung von
  Technik und Nachhaltigkeit}.
\newblock (German) conference talk at 14. VDV-Elektrobuskonferenz \& Fachmesse.

\bibitem[Erdeli{\'c} and Cari{\'c}, 2019]{EC19}
Erdeli{\'c}, T. and Cari{\'c}, T. (2019).
\newblock {A Survey on the Electric Vehicle Routing Problem: Variants and
  Solution Approaches}.
\newblock {\em Journal of Advanced Transportation}, 2019:1--48.

\bibitem[Froger et~al., 2022]{FJML22}
Froger, A., Jabali, O., Mendoza, J.~E., and Laporte, G. (2022).
\newblock {The Electric Vehicle Routing Problem with Capacitated Charging
  Stations}.
\newblock {\em Transportation Science}, 56(2):460--482.

\bibitem[Froger et~al., 2019]{FMJL19}
Froger, A., Mendoza, J.~E., Jabali, O., and Laporte, G. (2019).
\newblock Improved formulations and algorithmic components for the electric
  vehicle routing problem with nonlinear charging functions.
\newblock {\em Computers \& Operations Research}, 104:256--294.

\bibitem[{Gurobi Optimization, LLC}, 2024]{gurobi}
{Gurobi Optimization, LLC} (2024).
\newblock {Gurobi Optimizer Reference Manual}.

\bibitem[Hevia-Koch et~al., 2023]{IEA23}
Hevia-Koch, P., Wanner, B., Kuwahata, R., et~al. (2023).
\newblock {Electricity Grids and Secure Energy Transitions}.
\newblock Technical report, International Energy Agency, Paris.
\newblock last accessed June 2024.

\bibitem[Jahic et~al., 2019]{JES19}
Jahic, A., Eskander, M., and Schulz, D. (2019).
\newblock {Charging Schedule for Load Peak Minimization on Large-Scale Electric
  Bus Depots}.
\newblock {\em Applied Sciences}, 9(9):1748.

\bibitem[Jefferies and G{\"o}hlich, 2020]{JG20}
Jefferies, D. and G{\"o}hlich, D. (2020).
\newblock {A Comprehensive TCO Evaluation Method for Electric Bus Systems Based
  on Discrete-Event Simulation Including Bus Scheduling and Charging
  Infrastructure Optimisation}.
\newblock {\em World Electric Vehicle Journal}, 11.

\bibitem[Kim and Chung, 2023]{KC23}
Kim, Y.~J. and Chung, B.~D. (2023).
\newblock Energy consumption optimization for the electric vehicle routing
  problem with state-of-charge-dependent discharging rates.
\newblock {\em Journal of Cleaner Production}, 385:135703.

\bibitem[Kullman et~al., 2021]{KGM21}
Kullman, N.~D., Goodson, J.~C., and Mendoza, J.~E. (2021).
\newblock {Electric Vehicle Routing with Public Charging Stations}.
\newblock {\em Transportation Science}, 55(3):637--659.

\bibitem[Lee, 2021]{Lee21}
Lee, C. (2021).
\newblock An exact algorithm for the electric-vehicle routing problem with
  nonlinear charging time.
\newblock {\em Journal of the Operational Research Society}, 72(7):1461--1485.

\bibitem[Leptien, 2024]{Leptien24}
Leptien, A. (2024).
\newblock {Von Plug and Charge zum bidirektionalen Lade- und
  Energiemanagementsystem}.
\newblock {\em Der Nahverkehr Elektrobus-Spezial 2024}, 42:21–25.
\newblock (German).

\bibitem[Li et~al., 2019]{LLX19}
Li, L., Lo, H.~K., and Xiao, F. (2019).
\newblock Mixed bus fleet scheduling under range and refueling constraints.
\newblock {\em Transportation Research Part C: Emerging Technologies},
  104:443--462.

\bibitem[Liang et~al., 2021]{LDL21}
Liang, Y., Dabia, S., and Luo, Z. (2021).
\newblock {The Electric Vehicle Routing Problem with Nonlinear Charging
  Functions}.

\bibitem[L\"{o}bel et~al., 2023]{LBWatmos23}
L\"{o}bel, F., Bornd\"{o}rfer, R., and Weider, S. (2023).
\newblock {Non-Linear Charge Functions for Electric Vehicle Scheduling with
  Dynamic Recharge Rates}.
\newblock In Frigioni, D. and Schiewe, P., editors, {\em 23rd Symposium on
  Algorithmic Approaches for Transportation Modelling, Optimization, and
  Systems (ATMOS 2023)}, volume 115 of {\em Open Access Series in Informatics
  (OASIcs)}, pages 15:1--15:6, Dagstuhl, Germany. Schloss Dagstuhl --
  Leibniz-Zentrum f{\"u}r Informatik.

\bibitem[L\"{o}bel et~al., 2024]{LBWor24}
L\"{o}bel, F., Bornd\"{o}rfer, R., and Weider, S. (2024).
\newblock {Non-linear Battery Behavior in Electric Vehicle Scheduling
  Problems}.
\newblock In Voigt, G., Fliedner, M., Haase, K., Br{\"u}ggemann, W., Hoberg,
  K., and Meissner, J., editors, {\em Operations Research Proceedings 2023},
  Lecture Notes in Operations Research. Springer Cham.

\bibitem[L{\"o}sel, 2023]{Losel23}
L{\"o}sel, S. (2023).
\newblock {Elektrobusse im l{\"a}ndlichen Raum: VLP startet Sektorenkopplung
  Energiewirtschaft und Verkehr}.
\newblock {\em Der Nahverkehr Elektrobus-Spezial 2023}, 41:40–43.
\newblock (German).

\bibitem[Messaoudi and Oulamara, 2019]{MO19}
Messaoudi, B. and Oulamara, A. (2019).
\newblock {Electric Bus Scheduling and Optimal Charging}.
\newblock In Paternina-Arboleda, C. and Vo{\ss}, S., editors, {\em
  Computational Logistics}, pages 233--247. Springer International Publishing.

\bibitem[Montoya et~al., 2017]{MGMV17}
Montoya, A., Gu{\'e}ret, C., Mendoza, J.~E., and Villegas, J.~G. (2017).
\newblock The electric vehicle routing problem with nonlinear charging
  function.
\newblock {\em Transportation Research Part B: Methodological}, 103:87--110.
\newblock Green Urban Transportation.

\bibitem[Nei{\ss}endorfer, 2023]{Neissen23}
Nei{\ss}endorfer, M. (2023).
\newblock {Kiepe Electric l{\"a}dt E-Busse mit Bahnstrom besonders schnell und
  effizient}.
\newblock {\em Nahverkehrs-praxis}, 6-2023:60--61.
\newblock (German).

\bibitem[Olsen and Kliewer, 2020]{OK20}
Olsen, N. and Kliewer, N. (2020).
\newblock {Scheduling Electric Buses in Public Transport: Modeling of the
  Charging Process and Analysis of Assumptions}.
\newblock {\em Logistics Research}, 13(4).

\bibitem[Pelletier et~al., 2017]{PJLV17}
Pelletier, S., Jabali, O., Laporte, G., and Veneroni, M. (2017).
\newblock {Battery degradation and behaviour for electric vehicles: Review and
  numerical analyses of several models}.
\newblock {\em Transportation Research Part B: Methodological}, 103:158--187.

\bibitem[Perumal et~al., 2022]{PLL22}
Perumal, S.~S., Lusby, R.~M., and Larsen, J. (2022).
\newblock Electric bus planning \& scheduling: A review of related problems and
  methodologies.
\newblock {\em European Journal of Operational Research}.

\bibitem[Richarz, 2024]{Richarz24}
Richarz, M. (2024).
\newblock {E-Mobilit{\"a}t ist mehr als Busse - Ein Blick auf die
  Betriebsh{\"o}fe der Rheinbahn}.
\newblock {\em Der Nahverkehr Elektrobus-Spezial 2024}, 42:6–10.
\newblock (German).

\bibitem[Rothhardt and Suerkemper, 2023]{RS23}
Rothhardt, S. and Suerkemper, L. (2023).
\newblock {Europas modernster Betriebshof f{\"u}r Elektrobusse erfolgreich in
  Betrieb}.
\newblock {\em Nahverkehrs-praxis}, 6-2023:62--64.
\newblock (German).

\bibitem[Schulze and Lackner, 2023]{SL23}
Schulze, S. and Lackner, J. (2023).
\newblock {Mehr saubere Busse im {\"O}PNV}.
\newblock {\em Nahverkehrs-praxis}, 6-2023:50--52.
\newblock (German).

\bibitem[Sclar et~al., 2019]{SGCL19}
Sclar, R., Gorguinpour, C., Castellanos, S., and Li, X. (2019).
\newblock {Barriers to Adopting Electric Buses}.
\newblock Technical report, World Resource Institute, 10 g street ne, suite
  800, Washington, DC 20002, USA.
\newblock last accessed June 2024.

\bibitem[Sievers and Laumen, 2023]{Nuremberg23}
Sievers, M. and Laumen, A. (2023).
\newblock {Gro{\ss}es Potenzial f{\"u}r den Einsatz von eBussen}.
\newblock {\em Nahverkehrs-praxis}, 6-2023:42--73.
\newblock (German).

\bibitem[van Kooten~Niekerk et~al., 2017]{KAH17}
van Kooten~Niekerk, M.~E., van~den Akker, J.~M., and Hoogeveen, J.~A. (2017).
\newblock {Scheduling electric vehicles}.
\newblock {\em Public Transport}, 9:155--176.

\bibitem[Wu et~al., 2022]{WLLJ22}
Wu, W., Lin, Y., Liu, R., and Jin, W. (2022).
\newblock The multi-depot electric vehicle scheduling problem with power grid
  characteristics.
\newblock {\em Transportation Research Part B: Methodological}, 155:322--347.

\bibitem[Zhang et~al., 2021a]{ZLTDCGL21}
Zhang, A., Li, T., Tu, R., Dong, C., Chen, H., Gao, J., and Liu, Y. (2021a).
\newblock {The Effect of Nonlinear Charging Function and Line Change
  Constraints on Electric Bus Scheduling}.
\newblock {\em Promet - Traffic \& Transportation}, 33(4):527--538.

\bibitem[Zhang et~al., 2021b]{ZWQ21}
Zhang, L., Wang, S., and Qu, X. (2021b).
\newblock Optimal electric bus fleet scheduling considering battery degradation
  and non-linear charging profile.
\newblock {\em Transportation Research Part E: Logistics and Transportation
  Review}, 154:102445.

\bibitem[Zhou et~al., 2022]{ZMO22}
Zhou, Y., Meng, Q., and Ong, G.~P. (2022).
\newblock {Electric Bus Charging Scheduling for a Single Public Transport Route
  Considering Nonlinear Charging Profile and Battery Degradation Effect}.
\newblock {\em Transportation Research Part B: Methodological}, 159:49--75.

\end{thebibliography}

\section*{Appendix: Computational Results}

\begin{table}[H]
    \footnotesize
    \begin{tabularx}{\textwidth}{c X r r | c r r r r r r r r}
        \toprule
        name & start sol & $\numLinSegments$ & $\timestep$ & fs.? & $\#\course$ & obj. & bound & \makecell[c]{gap\\(\%)} & $\#\course$ & obj. & bound & \makecell[c]{gap\\(\%)} \\
        \midrule
        \multirow{24}{*}{A} & \multirow{12}{*}{-} & 2 & 60 & - & 7 & 5\,812.24 & 5\,239.18 & 9.86 & 7 & 5\,803.77 & 5\,241.23 & 9.69 \\
        & & 3 & 60 & - & 7 & 5\,805.48 & 5\,230.58 & 9.90 & 7 & 5\,803.77 & 5\,233.47 & 9.83 \\
        & & 4 & 60 & - & 7 & 5\,809.65 & 5\,233.54 & 9.92 & 7 & 5\,803.77 & 5\,236.30 & 9.78 \\
        & & 10 & 60 & - & 7 & 5\,799.53 & 5\,223.49 & 9.93 & 7 & 5\,799.37 & 5\,226.62 & 9.88 \\
        & & 2 & 300 & - & 7 & 5\,807.39 & 5\,262.62 & 9.38 & 7 & 5\,803.77 & 5\,266.47 & 9.26 \\
        & & 3 & 300 & - & 7 & 5\,809.95 & 5\,266.89 & 9.35 & 7 & 5\,808.17 & 5\,272.94 & 9.22 \\
        & & 4 & 300 & - & 7 & 5\,805.01 & 5\,257.54 & 9.43 & 7 & 5\,803.78 & 5\,260.57 & 9.36 \\
        & & 10 & 300 & - & 7 & 5\,805.03 & 5\,259.01 & 9.41 & 7 & 5\,803.77 & 5\,264.41 & 9.29 \\
        & & 2 & 600 & - & 7 & 5\,804.09 & 5\,279.96 & 9.03 & 7 & 5\,803.78 & 5\,292.58 & 8.81 \\
        & & 3 & 600 & - & 7 & 5\,803.77 & 5\,276.73 & 9.08 & 7 & 5\,803.77 & 5\,282.98 & 8.97 \\
        & & 4 & 600 & - & 7 & 5\,803.77 & 5\,273.79 & 9.13 & 7 & 5\,803.77 & 5\,280.11 & 9.02 \\
        & & 10 & 600 & - & 7 & 5\,807.53 & 5\,269.96 & 9.26 & 7 & 5\,803.78 & 5\,273.41 & 9.14 \\
        \cmidrule{2-13}
        & \multirow{12}{*}{5\,234.57} & 2 & 60 & $\times$ & 7 & 5\,812.24 & 5\,239.19 & 9.86 & 7 & 5\,803.77 & 5\,241.24 & 9.69 \\
        & & 3 & 60 & $\times$ & 7 & 5\,807.54 & 5\,235.88 & 9.84 & 7 & 5\,803.77 & 5\,238.88 & 9.73 \\
        & & 4 & 60 & $\times$ & 7 & 5\,809.65 & 5\,233.54 & 9.92 & 7 & 5\,803.77 & 5\,236.24 & 9.78 \\
        & & 10 & 60 & $\times$ & 7 & 5\,815.91 & 5\,237.15 & 9.95 & 7 & 5\,808.18 & 5\,240.38 & 9.78 \\
        & & 2 & 300 & $\times$ & 7 & 5\,807.39 & 5\,262.67 & 9.38 & 7 & 5\,803.77 & 5\,266.46 & 9.26 \\
        & & 3 & 300 & $\times$ & 7 & 5\,804.75 & 5\,262.90 & 9.33 & 7 & 5\,803.78 & 5\,269.71 & 9.20 \\
        & & 4 & 300 & $\times$ & 7 & 5\,806.41 & 5\,257.77 & 9.45 & 7 & 5\,803.77 & 5\,262.25 & 9.33 \\
        & & 10 & 300 & $\times$ & 7 & 5\,806.11 & 5\,259.01 & 9.42 & 7 & 5\,803.77 & 5\,264.40 & 9.29 \\
        & & 2 & 600 & $\times$ & 7 & 5\,809.73 & 5\,280.84 & 9.10 & 7 & 5\,803.77 & 5\,296.12 & 8.75 \\
        & & 3 & 600 & $\times$ & 7 & 5\,803.77 & 5\,276.73 & 9.08 & 7 & 5\,803.77 & 5\,282.98 & 8.97 \\
        & & 4 & 600 & $\times$ & 7 & 5\,803.77 & 5\,273.79 & 9.13 & 7 & 5\,803.77 & 5\,280.11 & 9.02 \\
        & & 10 & 600 & $\times$ & 7 & 5\,802.23 & 5\,266.48 & 9.23 & 7 & 5\,799.38 & 5\,270.54 & 9.12 \\
        \bottomrule
    \end{tabularx}
    \caption{Computational results testing the discretization parameter $\numLinSegments$ and $\timestep$.
    Each configuration was used on each instance once without and once with a start solution which is energy-feasible (and for A, C, G, J, L, O, and P (almost) optimal) under a linear charge curve approximation.
    Its objective is given in the third column.
    The fifth column indicates whether the given initial vehicle schedule was energy feasible under the corresponding configuration.
    The remaining columns give the fleet size ($\#\course$), objective value, lower bound, and integrality gap after 1 hour and after 12 hours runtime.}
    \label{tab:discParam1}
\end{table}

\begin{table}[H]
    \footnotesize
    \begin{tabularx}{\textwidth}{c X r r | c r r r r r r r r}
        \toprule
        name & start sol & $\numLinSegments$ & $\timestep$ & fs.? & $\#\course$ & obj. & bound & \makecell[c]{gap\\(\%)} & $\#\course$ & obj. & bound & \makecell[c]{gap\\(\%)} \\
        \midrule
        \multirow{24}{*}{B} & \multirow{12}{*}{-} & 2 & 60 & - & 8 & 6\,385.51 & 5\,804.78 & 9.09 & 8 & 6\,384.75 & 5\,839.84 & 8.53 \\
        & & 3 & 60 & - & 8 & 6\,384.75 & 5\,794.61 & 9.24 & 8 & 6\,384.75 & 5\,833.99 & 8.63 \\
        & & 4 & 60 & - & 8 & 6\,384.75 & 5\,801.63 & 9.13 & 8 & 6\,384.75 & 5\,831.44 & 8.67 \\
        & & 10 & 60 & - & 8 & 6\,400.60 & 5\,794.27 & 9.47 & 8 & 6\,384.75 & 5\,828.97 & 8.70 \\
        & & 2 & 300 & - & 8 & 6\,384.74 & 5\,838.09 & 8.56 & 8 & 6\,384.74 & 5\,871.09 & 8.04 \\
        & & 3 & 300 & - & 8 & 6\,386.61 & 5\,838.06 & 8.59 & 8 & 6\,384.75 & 5\,886.16 & 7.81 \\
        & & 4 & 300 & - & 8 & 6\,384.75 & 5\,833.40 & 8.64 & 8 & 6\,384.75 & 5\,875.80 & 7.97 \\
        & & 10 & 300 & - & 8 & 6\,384.91 & 5\,830.03 & 8.69 & 8 & 6\,384.75 & 5\,871.98 & 8.03 \\
        & & 2 & 600 & - & 8 & 6\,384.75 & 5\,857.91 & 8.25 & 8 & 6\,384.75 & 6\,202.41 & 2.86 \\
        & & 3 & 600 & - & 8 & 6\,384.75 & 5\,849.42 & 8.38 & 8 & 6\,384.75 & 5\,915.97 & 7.34 \\
        & & 4 & 600 & - & 8 & 6\,384.93 & 5\,855.97 & 8.28 & 8 & 6\,384.75 & 5\,924.71 & 7.21 \\
        & & 10 & 600 & - & 8 & 6\,384.92 & 5\,844.99 & 8.46 & 8 & 6\,384.75 & 5\,884.60 & 7.83 \\
        \cmidrule{2-13}
        & \multirow{12}{*}{6\,384.75} & 2 & 60 & $\times$ & 8 & 6\,397.39 & 5\,801.91 & 9.31 & 8 & 6\,384.75 & 5\,839.23 & 8.54 \\
        & & 3 & 60 & $\times$ & 8 & 6\,389.05 & 5\,806.71 & 9.11 & 8 & 6\,384.75 & 5\,840.79 & 8.52 \\
        & & 4 & 60 & $\times$ & 8 & 6\,397.25 & 5\,802.31 & 9.30 & 8 & 6\,384.75 & 5\,831.37 & 8.67 \\
        & & 10 & 60 & $\times$ & 8 & 6\,398.46 & 5\,788.20 & 9.54 & 8 & 6\,384.75 & 5\,818.85 & 8.86 \\
        & & 2 & 300 & $\times$ & 8 & 6\,384.74 & 5\,847.40 & 8.42 & 8 & 6\,384.74 & 6\,159.19 & 3.53 \\
        & & 3 & 300 & $\times$ & 8 & 6\,384.75 & 5\,839.09 & 8.55 & 8 & 6\,384.75 & 5\,882.80 & 7.86 \\
        & & 4 & 300 & $\times$ & 8 & 6\,384.75 & 5\,834.72 & 8.61 & 8 & 6\,384.75 & 5\,897.64 & 7.63 \\
        & & 10 & 300 & $\times$ & 8 & 6\,384.76 & 5\,821.89 & 8.82 & 8 & 6\,384.75 & 5\,872.06 & 8.03 \\
        & & 2 & 600 & $\times$ & 8 & 6\,386.06 & 5\,853.59 & 8.34 & 8 & 6\,384.75 & 6\,011.92 & 5.84 \\
        & & 3 & 600 & $\times$ & 8 & 6\,386.03 & 5\,851.88 & 8.36 & 8 & 6\,384.75 & 5\,889.11 & 7.76 \\
        & & 4 & 600 & $\times$ & 8 & 6\,384.75 & 5\,849.12 & 8.39 & 8 & 6\,384.75 & 6\,083.09 & 4.72 \\
        & & 10 & 600 & $\times$ & 8 & 6\,384.75 & 5\,839.54 & 8.54 & 8 & 6\,384.75 & 5\,883.15 & 7.86 \\
        \midrule
        \multirow{24}{*}{C} & \multirow{12}{*}{-} & 2 & 60 & - & 9 & 6\,125.52 & 6\,099.34 & 0.43 & 9 & 6\,122.82 & 6\,102.21 & 0.34 \\
        & & 3 & 60 & - & - & - & 6\,097.48 & - & 9 & 6\,119.87 & 6\,100.48 & 0.32 \\
        & & 4 & 60 & - & - & - & 6\,100.81 & - & 9 & 6\,124.83 & 6\,103.02 & 0.36 \\
        & & 10 & 60 & - & - & - & 6\,099.48 & - & 9 & 6\,123.13 & 6\,100.95 & 0.36 \\
        & & 2 & 300 & - & 9 & 6\,122.00 & 6\,098.57 & 0.38 & 9 & 6\,121.73 & 6\,104.82 & 0.28 \\
        & & 3 & 300 & - & 9 & 6\,125.36 & 6\,102.30 & 0.38 & 9 & 6\,125.00 & 6\,105.79 & 0.31 \\
        & & 4 & 300 & - & 9 & 6\,119.77 & 6\,097.89 & 0.36 & 9 & 6\,119.75 & 6\,103.70 & 0.26 \\
        & & 10 & 300 & - & 9 & 6\,121.18 & 6\,096.94 & 0.40 & 9 & 6\,119.40 & 6\,100.69 & 0.31 \\
        & & 2 & 600 & - & 9 & 6\,141.49 & 6\,096.96 & 0.73 & 9 & 6\,121.93 & 6\,101.76 & 0.33 \\
        & & 3 & 600 & - & 9 & 6\,122.32 & 6\,099.66 & 0.37 & 9 & 6\,121.70 & 6\,108.63 & 0.21 \\
        & & 4 & 600 & - & 9 & 6\,122.07 & 6\,099.10 & 0.38 & 9 & 6\,121.82 & 6\,107.27 & 0.24 \\
        & & 10 & 600 & - & 10 & 6\,462.61 & 6\,096.47 & 5.67 & 9 & 6\,120.45 & 6\,100.92 & 0.32 \\
        \cmidrule{2-13}
        & \multirow{12}{*}{6\,124.55} & 2 & 60 & $\times$ & - & - & 6\,098.40 & - & 9 & 6\,121.98 & 6\,100.77 & 0.35 \\
        & & 3 & 60 & $\times$ & 9 & 6\,123.73 & 6\,101.12 & 0.37 & 9 & 6\,122.32 & 6\,103.02 & 0.32 \\
        & & 4 & 60 & $\times$ & 9 & 6\,123.23 & 6\,098.94 & 0.40 & 9 & 6\,121.74 & 6\,100.54 & 0.35 \\
        & & 10 & 60 & $\times$ & - & - & 6\,099.47 & - & 9 & 6\,123.13 & 6\,100.99 & 0.36 \\
        & & 2 & 300 & $\times$ & 9 & 6\,122.75 & 6\,098.30 & 0.40 & 9 & 6\,122.03 & 6\,103.62 & 0.30 \\
        & & 3 & 300 & $\times$ & 9 & 6\,126.53 & 6\,100.37 & 0.43 & 9 & 6\,123.12 & 6\,106.33 & 0.27 \\
        & & 4 & 300 & $\times$ & 9 & 6\,120.68 & 6\,096.86 & 0.39 & 9 & 6\,119.82 & 6\,098.46 & 0.35 \\
        & & 10 & 300 & $\times$ & 9 & 6\,124.48 & 6\,099.13 & 0.41 & 9 & 6\,122.12 & 6\,103.37 & 0.31 \\
        & & 2 & 600 & $\times$ & 9 & 6\,125.07 & 6\,099.54 & 0.42 & 9 & 6\,123.29 & 6\,106.06 & 0.28 \\
        & & 3 & 600 & $\times$ & 9 & 6\,121.04 & 6\,097.49 & 0.38 & 9 & 6\,119.90 & 6\,103.01 & 0.28 \\
        & & 4 & 600 & $\times$ & 9 & 6\,121.25 & 6\,095.79 & 0.42 & 9 & 6\,120.95 & 6\,104.42 & 0.27 \\
        & & 10 & 600 & $\times$ & 9 & 6\,125.99 & 6\,102.09 & 0.39 & 9 & 6\,124.66 & 6\,110.15 & 0.24 \\
        \bottomrule
    \end{tabularx}
    \caption{\Cref{tab:discParam1} continued.}
    \label{tab:discParam2}
\end{table}

\begin{table}[H]
    \scriptsize
    \begin{tabularx}{\textwidth}{c X r r | c r r r r r r r r}
        \toprule
        name & start sol & $\numLinSegments$ & $\timestep$ & fs.? & $\#\course$ & obj. & bound & \makecell[c]{gap\\(\%)} & $\#\course$ & obj. & bound & \makecell[c]{gap\\(\%)} \\
        \midrule
        \multirow{24}{*}{D} & \multirow{12}{*}{-} & 2 & 60 & - & 9 & 6\,803.34 & 6\,258.32 & 8.01 & 9 & 6\,803.34 & 6\,338.21 & 6.84 \\
        & & 3 & 60 & - & 9 & 6\,803.35 & 6\,257.48 & 8.02 & 9 & 6\,803.35 & 6\,271.98 & 7.81 \\
        & & 4 & 60 & - & 9 & 6\,803.34 & 6\,252.65 & 8.09 & 9 & 6\,803.34 & 6\,281.75 & 7.67 \\
        & & 10 & 60 & - & 9 & 6\,817.14 & 6\,254.44 & 8.25 & 9 & 6\,803.37 & 6\,266.40 & 7.89 \\
        & & 2 & 300 & - & 9 & 6\,803.39 & 6\,283.30 & 7.64 & 9 & 6\,803.34 & 6\,778.61 & 0.36 \\
        & & 3 & 300 & - & 9 & 6\,803.46 & 6\,278.81 & 7.71 & 9 & 6\,803.34 & 6\,483.49 & 4.70 \\
        & & 4 & 300 & - & 9 & 6\,803.34 & 6\,273.80 & 7.78 & 9 & 6\,803.34 & 6\,416.61 & 5.68 \\
        & & 10 & 300 & - & 9 & 6\,803.34 & 6\,274.04 & 7.78 & 9 & 6\,803.34 & 6\,777.99 & 0.37 \\
        & & 2 & 600 & - & 9 & 6\,804.59 & 6\,304.78 & 7.35 & 9 & 6\,803.34 & 6\,412.40 & 5.75 \\
        & & 3 & 600 & - & 9 & 6\,803.34 & 6\,296.46 & 7.45 & 9 & 6\,803.34 & 6\,385.30 & 6.14 \\
        & & 4 & 600 & - & 9 & 6\,803.35 & 6\,284.30 & 7.63 & 9 & 6\,803.35 & 6\,779.20 & 0.35 \\
        & & 10 & 600 & - & 9 & 6\,803.35 & 6\,305.25 & 7.32 & 9 & 6\,803.35 & 6\,497.18 & 4.50 \\
        \cmidrule{2-13}
        & \multirow{12}{*}{6\,803.34} & 2 & 60 & $\checkmark$ & 9 & 6\,803.34 & 6\,262.34 & 7.95 & 9 & 6\,803.34 & 6\,293.70 & 7.49 \\
        & & 3 & 60 & $\checkmark$ & 9 & 6\,803.34 & 6\,256.28 & 8.04 & 9 & 6\,803.34 & 6\,273.79 & 7.78 \\
        & & 4 & 60 & $\checkmark$ & 9 & 6\,803.34 & 6\,252.17 & 8.10 & 9 & 6\,803.34 & 6\,265.57 & 7.90 \\
        & & 10 & 60 & $\checkmark$ & 9 & 6\,803.34 & 6\,253.83 & 8.08 & 9 & 6\,803.34 & 6\,259.71 & 7.99 \\
        & & 2 & 300 & $\checkmark$ & 9 & 6\,803.34 & 6\,274.84 & 7.77 & 9 & 6\,803.34 & 6\,777.85 & 0.37 \\
        & & 3 & 300 & $\checkmark$ & 9 & 6\,803.34 & 6\,276.83 & 7.74 & 9 & 6\,803.34 & 6\,492.28 & 4.57 \\
        & & 4 & 300 & $\checkmark$ & 9 & 6\,803.34 & 6\,274.51 & 7.77 & 9 & 6\,803.34 & 6\,384.46 & 6.16 \\
        & & 10 & 300 & $\checkmark$ & 9 & 6\,803.34 & 6\,273.81 & 7.78 & 9 & 6\,803.34 & 6\,620.50 & 2.69 \\
        & & 2 & 600 & $\checkmark$ & 9 & 6\,803.34 & 6\,309.95 & 7.25 & 9 & 6\,803.34 & 6\,779.09 & 0.36 \\
        & & 3 & 600 & $\checkmark$ & 9 & 6\,803.34 & 6\,323.76 & 7.05 & 9 & 6\,803.34 & 6\,778.88 & 0.36 \\
        & & 4 & 600 & $\checkmark$ & 9 & 6\,803.34 & 6\,308.10 & 7.28 & 9 & 6\,803.34 & 6\,779.78 & 0.35 \\
        & & 10 & 600 & $\checkmark$ & 9 & 6\,803.34 & 6\,301.90 & 7.37 & 9 & 6\,803.34 & 6\,778.91 & 0.36 \\
        \midrule
        \multirow{24}{*}{E} & \multirow{12}{*}{-} & 2 & 60 & - & 16 & 12\,561.75 & 10\,073.55 & 19.81 & 13 & 10\,742.05 & 10\,097.55 & 6.00 \\
        & & 3 & 60 & - & 16 & 12\,542.99 & 10\,069.69 & 19.72 & 13 & 10\,745.59 & 10\,083.09 & 6.17 \\
        & & 4 & 60 & - & 14 & 11\,387.94 & 10\,074.84 & 11.53 & 14 & 11\,326.24 & 10\,097.24 & 10.85 \\
        & & 10 & 60 & - & 16 & 12\,559.85 & 10\,069.85 & 19.83 & 14 & 11\,350.05 & 10\,083.25 & 11.16 \\
        & & 2 & 300 & - & 14 & 11\,309.28 & 10\,139.58 & 10.34 & 14 & 11\,299.08 & 10\,150.68 & 10.16 \\
        & & 3 & 300 & - & 14 & 11\,335.98 & 10\,101.78 & 10.89 & 13 & 10\,732.38 & 10\,141.88 & 5.50 \\
        & & 4 & 300 & - & 14 & 11\,323.44 & 10\,115.64 & 10.67 & 13 & 10\,749.94 & 10\,137.54 & 5.70 \\
        & & 10 & 300 & - & 14 & 11\,402.94 & 10\,088.14 & 11.53 & 13 & 10\,741.44 & 10\,137.94 & 5.62 \\
        & & 2 & 600 & - & 14 & 11\,311.03 & 10\,172.23 & 10.07 & 14 & 11\,293.23 & 10\,199.93 & 9.68 \\
        & & 3 & 600 & - & 13 & 10\,774.30 & 10\,156.10 & 5.74 & 13 & 10\,732.40 & 10\,172.60 & 5.22 \\
        & & 4 & 600 & - & 13 & 10\,732.63 & 10\,145.13 & 5.47 & 13 & 10\,732.43 & 10\,164.63 & 5.29 \\
        & & 10 & 600 & - & 14 & 11\,318.88 & 10\,150.58 & 10.32 & 13 & 10\,794.38 & 10\,170.28 & 5.78 \\
        \cmidrule{2-13}
        & \multirow{12}{*}{10\,707.00} & 2 & 60 & $\times$ & 14 & 11\,383.73 & 10\,074.53 & 11.50 & 13 & 10\,744.13 & 10\,101.53 & 5.98 \\
        & & 3 & 60 & $\times$ & 15 & 11\,949.95 & 10\,070.75 & 15.73 & 14 & 11\,309.75 & 10\,091.55 & 10.77 \\
        & & 4 & 60 & $\times$ & 16 & 12\,564.80 & 10\,072.20 & 19.84 & 14 & 11\,320.40 & 10\,085.30 & 10.91 \\
        & & 10 & 60 & $\times$ & 16 & 12\,573.45 & 10\,069.85 & 19.91 & 14 & 11\,350.05 & 10\,083.25 & 11.16 \\
        & & 2 & 300 & $\times$ & 14 & 11\,333.21 & 10\,099.91 & 10.88 & 13 & 10\,755.31 & 10\,145.81 & 5.67 \\
        & & 3 & 300 & $\times$ & 14 & 11\,342.54 & 10\,100.24 & 10.95 & 14 & 11\,293.24 & 10\,141.84 & 10.20 \\
        & & 4 & 300 & $\times$ & 14 & 11\,327.96 & 10\,089.66 & 10.93 & 13 & 10\,730.96 & 10\,144.26 & 5.47 \\
        & & 10 & 300 & $\times$ & 14 & 11\,344.46 & 10\,092.46 & 11.04 & 13 & 10\,742.06 & 10\,136.66 & 5.64 \\
        & & 2 & 600 & $\times$ & 14 & 11\,342.01 & 10\,171.71 & 10.32 & 14 & 11\,295.91 & 10\,202.21 & 9.68 \\
        & & 3 & 600 & $\times$ & 14 & 11\,338.67 & 10\,156.17 & 10.43 & 14 & 11\,317.77 & 10\,187.57 & 9.99 \\
        & & 4 & 600 & $\times$ & 14 & 11\,311.53 & 10\,150.83 & 10.26 & 13 & 10\,732.43 & 10\,169.63 & 5.24 \\
        & & 10 & 600 & $\times$ & 14 & 11\,315.58 & 10\,145.28 & 10.34 & 13 & 10\,724.48 & 10\,167.78 & 5.19 \\
        \bottomrule
    \end{tabularx}
    \caption{\Cref{tab:discParam2} continued.}
    \label{tab:discParam3}
\end{table}

\begin{table}[H]
    \footnotesize
    \begin{tabularx}{\textwidth}{c X r r | c r r r r r r r r}
        \toprule
        name & start sol & $\numLinSegments$ & $\timestep$ & fs.? & $\#\course$ & obj. & bound & \makecell[c]{gap\\(\%)} & $\#\course$ & obj. & bound & \makecell[c]{gap\\(\%)} \\
        \midrule
        \multirow{24}{*}{F} & \multirow{12}{*}{-} & 2 & 60 & - & 10 & 8\,182.02 & 7\,311.95 & 10.63 & 10 & 8\,173.54 & 7\,385.21 & 9.64 \\
        & & 3 & 60 & - & 10 & 8\,186.55 & 7\,314.33 & 10.65 & 10 & 8\,181.23 & 7\,538.34 & 7.86 \\
        & & 4 & 60 & - & 10 & 8\,187.13 & 7\,308.45 & 10.73 & 10 & 8\,177.35 & 7\,331.99 & 10.34 \\
        & & 10 & 60 & - & 10 & 8\,190.12 & 7\,297.19 & 10.90 & 9 & 7\,604.27 & 7\,522.26 & 1.08 \\
        & & 2 & 300 & - & 10 & 8\,176.59 & 7\,337.10 & 10.27 & 10 & 8\,173.22 & 7\,538.35 & 7.77 \\
        & & 3 & 300 & - & 10 & 8\,176.35 & 7\,333.95 & 10.30 & 10 & 8\,169.66 & 7\,538.65 & 7.72 \\
        & & 4 & 300 & - & 10 & 8\,181.04 & 7\,327.61 & 10.43 & 10 & 8\,173.54 & 7\,503.49 & 8.20 \\
        & & 10 & 300 & - & 10 & 8\,185.56 & 7\,337.60 & 10.36 & 10 & 8\,181.22 & 7\,545.11 & 7.78 \\
        & & 2 & 600 & - & 10 & 8\,174.06 & 7\,383.27 & 9.67 & 10 & 8\,173.54 & 7\,544.44 & 7.70 \\
        & & 3 & 600 & - & 10 & 8\,181.23 & 7\,377.55 & 9.82 & 10 & 8\,181.23 & 7\,552.38 & 7.69 \\
        & & 4 & 600 & - & 10 & 8\,176.84 & 7\,370.96 & 9.86 & 10 & 8\,173.54 & 7\,546.50 & 7.67 \\
        & & 10 & 600 & - & 10 & 8\,170.22 & 7\,362.84 & 9.88 & 10 & 8\,169.67 & 7\,494.60 & 8.26 \\
        \cmidrule{2-13}
        & \multirow{12}{*}{8\,173.54} & 2 & 60 & $\times$ & 10 & 8\,187.21 & 7\,316.52 & 10.63 & 10 & 8\,181.23 & 7\,412.52 & 9.40 \\
        & & 3 & 60 & $\times$ & 10 & 8\,171.29 & 7\,310.34 & 10.54 & 10 & 8\,165.86 & 7\,518.07 & 7.93 \\
        & & 4 & 60 & $\times$ & 10 & 8\,173.77 & 7\,309.89 & 10.57 & 10 & 8\,173.21 & 7\,502.05 & 8.21 \\
        & & 10 & 60 & $\times$ & 10 & 8\,179.11 & 7\,303.50 & 10.71 & 10 & 8\,173.54 & 7\,323.87 & 10.40 \\
        & & 2 & 300 & $\times$ & 10 & 8\,169.98 & 7\,338.66 & 10.18 & 10 & 8\,165.86 & 7\,535.49 & 7.72 \\
        & & 3 & 300 & $\times$ & 10 & 8\,174.34 & 7\,340.90 & 10.20 & 10 & 8\,173.54 & 7\,540.89 & 7.74 \\
        & & 4 & 300 & $\times$ & 10 & 8\,174.61 & 7\,334.42 & 10.28 & 10 & 8\,173.55 & 7\,543.75 & 7.71 \\
        & & 10 & 300 & $\times$ & 10 & 8\,177.62 & 7\,322.31 & 10.46 & 10 & 8\,173.55 & 7\,376.18 & 9.76 \\
        & & 2 & 600 & $\times$ & 10 & 8\,177.38 & 7\,385.14 & 9.69 & 10 & 8\,173.55 & 7\,544.49 & 7.70 \\
        & & 3 & 600 & $\times$ & 10 & 8\,169.66 & 7\,373.09 & 9.75 & 10 & 8\,169.66 & 7\,399.29 & 9.43 \\
        & & 4 & 600 & $\times$ & 10 & 8\,173.90 & 7\,368.06 & 9.86 & 10 & 8\,173.54 & 7\,430.24 & 9.09 \\
        & & 10 & 600 & $\times$ & 10 & 8\,176.77 & 7\,359.60 & 9.99 & 10 & 8\,173.55 & 7\,543.92 & 7.70 \\
        \midrule
        \multirow{24}{*}{G} & \multirow{12}{*}{-} & 2 & 60 & - & 9 & 7\,504.67 & 7\,083.41 & 5.61 & 9 & 7\,504.67 & 7\,496.00 & 0.12 \\
        & & 3 & 60 & - & 9 & 7\,504.67 & 7\,067.70 & 5.82 & 9 & 7\,504.67 & 7\,498.28 & 0.09 \\
        & & 4 & 60 & - & 9 & 7\,504.67 & 7\,061.68 & 5.90 & 9 & 7\,504.67 & 7\,498.47 & 0.08 \\
        & & 10 & 60 & - & 9 & 7\,504.67 & 7\,059.07 & 5.94 & 9 & 7\,504.67 & 7\,113.09 & 5.22 \\
        & & 2 & 300 & - & 9 & 7\,504.67 & 7\,227.47 & 3.69 & 9 & 7\,504.67 & 7\,498.86 & 0.08 \\
        & & 3 & 300 & - & 9 & 7\,504.66 & 7\,495.59 & 0.12 & 9 & 7\,504.66 & 7\,499.00 & 0.08 \\
        & & 4 & 300 & - & 9 & 7\,504.67 & 7\,171.01 & 4.45 & 9 & 7\,504.67 & 7\,498.94 & 0.08 \\
        & & 10 & 300 & - & 9 & 7\,504.67 & 7\,172.35 & 4.43 & 9 & 7\,504.67 & 7\,498.78 & 0.08 \\
        & & 2 & 600 & - & 9 & 7\,504.67 & 7\,460.66 & 0.59 & 9 & 7\,504.67 & 7\,501.09 & 0.05 \\
        & & 3 & 600 & - & 9 & 7\,504.67 & 7\,463.28 & 0.55 & 9 & 7\,504.67 & 7\,500.60 & 0.05 \\
        & & 4 & 600 & - & 9 & 7\,504.93 & 7\,446.79 & 0.77 & 9 & 7\,504.67 & 7\,500.95 & 0.05 \\
        & & 10 & 600 & - & 9 & 7\,504.67 & 7\,463.47 & 0.55 & 9 & 7\,504.67 & 7\,500.92 & 0.05 \\
        \cmidrule{2-13}
        & \multirow{12}{*}{7\,504.67} & 2 & 60 & $\checkmark$ & 9 & 7\,504.67 & 7\,078.45 & 5.68 & 9 & 7\,504.67 & 7\,495.81 & 0.12 \\
        & & 3 & 60 & $\checkmark$ & 9 & 7\,504.67 & 7\,066.66 & 5.84 & 9 & 7\,504.67 & 7\,141.03 & 4.85 \\
        & & 4 & 60 & $\checkmark$ & 9 & 7\,504.67 & 7\,065.86 & 5.85 & 9 & 7\,504.67 & 7\,498.06 & 0.09 \\
        & & 10 & 60 & $\checkmark$ & 9 & 7\,504.67 & 7\,058.59 & 5.94 & 9 & 7\,504.67 & 7\,391.31 & 1.51 \\
        & & 2 & 300 & $\checkmark$ & 9 & 7\,504.67 & 7\,334.65 & 2.27 & 9 & 7\,504.67 & 7\,498.98 & 0.08 \\
        & & 3 & 300 & $\checkmark$ & 9 & 7\,504.66 & 7\,493.39 & 0.15 & 9 & 7\,504.66 & 7\,498.21 & 0.09 \\
        & & 4 & 300 & $\checkmark$ & 9 & 7\,504.67 & 7\,173.37 & 4.41 & 9 & 7\,504.67 & 7\,497.77 & 0.09 \\
        & & 10 & 300 & $\checkmark$ & 9 & 7\,504.67 & 7\,167.69 & 4.49 & 9 & 7\,504.67 & 7\,246.26 & 3.44 \\
        & & 2 & 600 & $\times$ & 9 & 7\,504.67 & 7\,460.95 & 0.58 & 9 & 7\,504.67 & 7\,501.37 & 0.04 \\
        & & 3 & 600 & $\times$ & 9 & 7\,504.67 & 7\,462.05 & 0.57 & 9 & 7\,504.67 & 7\,500.77 & 0.05 \\
        & & 4 & 600 & $\times$ & 9 & 7\,504.67 & 7\,446.09 & 0.78 & 9 & 7\,504.67 & 7\,495.48 & 0.12 \\
        & & 10 & 600 & $\times$ & 9 & 7\,504.67 & 7\,459.27 & 0.60 & 9 & 7\,504.67 & 7\,501.17 & 0.05 \\
        \bottomrule
    \end{tabularx}
    \caption{\Cref{tab:discParam3} continued.}
    \label{tab:discParam4}
\end{table}

\begin{table}[H]
    \footnotesize
    \begin{tabularx}{\textwidth}{c X r r | c r r r r r r r r}
        \toprule
        name & start sol & $\numLinSegments$ & $\timestep$ & fs.? & $\#\course$ & obj. & bound & \makecell[c]{gap\\(\%)} & $\#\course$ & obj. & bound & \makecell[c]{gap\\(\%)} \\
        \midrule
        \multirow{24}{*}{H} & \multirow{12}{*}{-} & 2 & 60 & - & - & - & 14\,449.53 & - & 18 & 14\,755.43 & 14\,452.13 & 2.06 \\
        & & 3 & 60 & - & - & - & 14\,434.41 & - & - & - & 14\,439.31 & - \\
        & & 4 & 60 & - & - & - & 14\,421.31 & - & - & - & 14\,439.03 & - \\
        & & 10 & 60 & - & - & - & 14\,421.20 & - & - & - & 14\,438.13 & - \\
        & & 2 & 300 & - & - & - & 14\,441.20 & - & - & - & 14\,456.76 & - \\
        & & 3 & 300 & - & - & - & 14\,440.56 & - & 18 & 14\,712.46 & 14\,454.96 & 1.75 \\
        & & 4 & 300 & - & - & - & 14\,436.31 & - & - & - & 14\,450.48 & - \\
        & & 10 & 300 & - & - & - & 14\,453.07 & - & 18 & 14\,993.37 & 14\,468.27 & 3.50 \\
        & & 2 & 600 & - & - & - & 14\,452.20 & - & - & - & 14\,468.08 & - \\
        & & 3 & 600 & - & - & - & 14\,450.60 & - & - & - & 14\,464.83 & - \\
        & & 4 & 600 & - & - & - & 14\,449.41 & - & - & - & 14\,462.96 & - \\
        & & 10 & 600 & - & - & - & 14\,445.91 & - & - & - & 14\,461.11 & - \\
        \cmidrule{2-13}
        & \multirow{12}{*}{14\,650.00} & 2 & 60 & $\times$ & - & - & 14\,423.31 & - & - & - & 14\,441.62 & - \\
        & & 3 & 60 & $\times$ & - & - & 14\,423.81 & - & - & - & 14\,440.63 & - \\
        & & 4 & 60 & $\times$ & - & - & 14\,421.31 & - & - & - & 14\,439.03 & - \\
        & & 10 & 60 & $\times$ & - & - & 14\,421.20 & - & - & - & 14\,438.13 & - \\
        & & 2 & 300 & $\times$ & - & - & 14\,443.64 & - & 18 & 14\,710.44 & 14\,455.34 & 1.73 \\
        & & 3 & 300 & $\times$ & - & - & 14\,438.20 & - & - & - & 14\,449.73 & - \\
        & & 4 & 300 & $\times$ & - & - & 14\,448.99 & - & 18 & 14\,719.39 & 14\,462.29 & 1.75 \\
        & & 10 & 300 & $\times$ & - & - & 14\,436.01 & - & - & - & 14\,449.31 & - \\
        & & 2 & 600 & $\times$ & - & - & 14\,453.81 & - & 18 & 14\,826.70 & 14\,469.51 & 2.41 \\
        & & 3 & 600 & $\times$ & - & - & 14\,450.20 & - & - & - & 14\,465.77 & - \\
        & & 4 & 600 & $\times$ & - & - & 14\,448.81 & - & - & - & 14\,461.93 & - \\
        & & 10 & 600 & $\times$ & - & - & 14\,445.20 & - & - & - & 14\,460.82 & - \\
        \midrule
        \multirow{24}{*}{I} & \multirow{12}{*}{-} & 2 & 60 & - & - & - & 13\,813.22 & - & 18 & 14\,724.32 & 13\,822.42 & 6.13 \\
        & & 3 & 60 & - & - & - & 13\,804.98 & - & 18 & 14\,602.48 & 13\,814.38 & 5.40 \\
        & & 4 & 60 & - & - & - & 13\,784.18 & - & - & - & 13\,794.05 & - \\
        & & 10 & 60 & - & - & - & 13\,782.28 & - & - & - & 13\,792.58 & - \\
        & & 2 & 300 & - & - & - & 13\,811.08 & - & 18 & 14\,627.68 & 13\,822.08 & 5.51 \\
        & & 3 & 300 & - & - & - & 13\,819.99 & - & 18 & 14\,644.99 & 13\,829.59 & 5.57 \\
        & & 4 & 300 & - & - & - & 13\,812.71 & - & 18 & 14\,599.91 & 13\,823.31 & 5.32 \\
        & & 10 & 300 & - & - & - & 13\,816.58 & - & 18 & 14\,669.78 & 13\,827.78 & 5.74 \\
        & & 2 & 600 & - & - & - & 13\,850.26 & - & 18 & 15\,140.76 & 13\,861.76 & 8.45 \\
        & & 3 & 600 & - & - & - & 13\,842.88 & - & - & - & 13\,855.97 & - \\
        & & 4 & 600 & - & - & - & 13\,854.55 & - & 18 & 14\,597.25 & 13\,868.15 & 4.99 \\
        & & 10 & 600 & - & - & - & 13\,853.48 & - & - & - & 13\,857.24 & - \\
        \cmidrule{2-13}
        & \multirow{12}{*}{13\,995.10} & 2 & 60 & $\times$ & - & - & 13\,812.94 & - & 18 & 14\,723.24 & 13\,822.14 & 6.12 \\
        & & 3 & 60 & $\times$ & - & - & 13\,804.98 & - & 18 & 14\,602.48 & 13\,814.38 & 5.40 \\
        & & 4 & 60 & $\times$ & - & - & 13\,788.47 & - & - & - & 13\,798.34 & - \\
        & & 10 & 60 & $\times$ & - & - & 13\,786.57 & - & - & - & 13\,796.87 & - \\
        & & 2 & 300 & $\times$ & - & - & 13\,825.08 & - & 18 & 14\,620.08 & 13\,835.78 & 5.36 \\
        & & 3 & 300 & $\times$ & - & - & 13\,813.93 & - & 18 & 14\,563.43 & 13\,824.33 & 5.08 \\
        & & 4 & 300 & $\times$ & - & - & 13\,818.43 & - & 18 & 14\,575.53 & 13\,828.63 & 5.12 \\
        & & 10 & 300 & $\times$ & - & - & 13\,817.61 & - & 18 & 14\,612.11 & 13\,828.61 & 5.36 \\
        & & 2 & 600 & $\times$ & - & - & 13\,848.37 & - & 18 & 14\,712.97 & 13\,863.47 & 5.77 \\
        & & 3 & 600 & $\times$ & - & - & 13\,860.28 & - & 18 & 15\,105.28 & 13\,874.08 & 8.15 \\
        & & 4 & 600 & $\times$ & - & - & 13\,864.36 & - & 18 & 14\,944.46 & 13\,876.16 & 7.15 \\
        & & 10 & 600 & $\times$ & - & - & 13\,860.70 & - & 18 & 14\,845.40 & 13\,875.20 & 6.54 \\
        \bottomrule
    \end{tabularx}
    \caption{\Cref{tab:discParam4} continued.}
    \label{tab:discParam5}
\end{table}

\begin{table}[H]
    \scriptsize
    \begin{tabularx}{\textwidth}{c X r r | c r r r r r r r r}
        \toprule
        name & start sol & $\numLinSegments$ & $\timestep$ & fs.? & $\#\course$ & obj. & bound & \makecell[c]{gap\\(\%)} & $\#\course$ & obj. & bound & \makecell[c]{gap\\(\%)} \\
        \midrule
        \multirow{24}{*}{J} & \multirow{12}{*}{-} & 2 & 60 & - & 42 & 124\,683.02 & 122\,211.02 & 1.98 & 37 & 122\,971.02 & 122\,211.02 & 0.62 \\
        & & 3 & 60 & - & 46 & 126\,996.67 & 122\,212.67 & 3.77 & 39 & 124\,257.67 & 122\,217.67 & 1.64 \\
        & & 4 & 60 & - & 43 & 125\,945.04 & 122\,195.04 & 2.98 & 38 & 123\,742.04 & 122\,203.04 & 1.24 \\
        & & 10 & 60 & - & 45 & 127\,122.05 & 122\,179.05 & 3.89 & 42 & 125\,713.05 & 122\,196.05 & 2.80 \\
        & & 2 & 300 & - & 38 & 123\,401.67 & 122\,219.67 & 0.96 & 37 & 122\,859.67 & 122\,219.67 & 0.52 \\
        & & 3 & 300 & - & 37 & 123\,609.93 & 122\,210.93 & 1.13 & 37 & 122\,995.93 & 122\,210.93 & 0.64 \\
        & & 4 & 300 & - & 41 & 124\,737.05 & 122\,215.05 & 2.02 & 37 & 122\,883.05 & 122\,307.05 & 0.47 \\
        & & 10 & 300 & - & 40 & 125\,075.44 & 122\,199.44 & 2.30 & 37 & 122\,944.44 & 122\,199.44 & 0.61 \\
        & & 2 & 600 & - & 37 & 123\,388.88 & 122\,223.88 & 0.94 & 37 & 122\,899.88 & 122\,224.88 & 0.55 \\
        & & 3 & 600 & - & 37 & 123\,794.53 & 122\,217.53 & 1.27 & 37 & 122\,832.53 & 122\,217.53 & 0.50 \\
        & & 4 & 600 & - & 38 & 123\,546.42 & 122\,206.42 & 1.08 & 37 & 122\,841.42 & 122\,206.42 & 0.52 \\
        & & 10 & 600 & - & 38 & 123\,952.29 & 122\,205.29 & 1.41 & 37 & 123\,021.29 & 122\,205.29 & 0.66 \\
        \cmidrule{2-13}
        & \multirow{12}{*}{122\,836.00} & 2 & 60 & $\checkmark$ & 37 & 122\,836.00 & 122\,201.98 & 0.52 & 37 & 122\,805.98 & 122\,201.98 & 0.49 \\
        & & 3 & 60 & $\checkmark$ & 37 & 122\,829.63 & 122\,222.63 & 0.49 & 37 & 122\,807.63 & 122\,224.63 & 0.47 \\
        & & 4 & 60 & $\checkmark$ & 37 & 122\,830.80 & 122\,187.80 & 0.52 & 37 & 122\,828.80 & 122\,193.80 & 0.52 \\
        & & 10 & 60 & $\checkmark$ & 37 & 122\,832.62 & 122\,180.62 & 0.53 & 37 & 122\,829.62 & 122\,197.62 & 0.51 \\
        & & 2 & 300 & $\times$ & 38 & 123\,403.11 & 122\,221.11 & 0.96 & 37 & 122\,873.11 & 122\,221.11 & 0.53 \\
        & & 3 & 300 & $\times$ & 37 & 123\,609.93 & 122\,210.93 & 1.13 & 37 & 122\,995.93 & 122\,210.93 & 0.64 \\
        & & 4 & 300 & $\times$ & 41 & 124\,737.05 & 122\,215.05 & 2.02 & 37 & 122\,883.05 & 122\,306.05 & 0.47 \\
        & & 10 & 300 & $\times$ & 40 & 125\,075.44 & 122\,199.44 & 2.30 & 37 & 122\,944.44 & 122\,199.44 & 0.61 \\
        & & 2 & 600 & $\times$ & 37 & 123\,388.88 & 122\,223.88 & 0.94 & 37 & 122\,899.88 & 122\,224.88 & 0.55 \\
        & & 3 & 600 & $\times$ & 37 & 123\,794.53 & 122\,217.53 & 1.27 & 37 & 122\,832.53 & 122\,217.53 & 0.50 \\
        & & 4 & 600 & $\times$ & 38 & 123\,546.42 & 122\,206.42 & 1.08 & 37 & 122\,841.42 & 122\,206.42 & 0.52 \\
        & & 10 & 600 & $\times$ & 39 & 124\,252.27 & 122\,205.27 & 1.65 & 37 & 123\,021.27 & 122\,205.27 & 0.66 \\
        \midrule
        \multirow{24}{*}{K} & \multirow{12}{*}{-} & 2 & 60 & - & - & - & 64\,034.32 & - & 188 & 134\,181.92 & 64\,287.62 & 52.09 \\
        & & 3 & 60 & - & - & - & 64\,036.96 & - & 235 & 163\,065.16 & 64\,280.76 & 60.58 \\
        & & 4 & 60 & - & - & - & 64\,023.41 & - & 236 & 163\,825.41 & 64\,255.61 & 60.78 \\
        & & 10 & 60 & - & - & - & 64\,003.62 & - & 241 & 166\,908.22 & 64\,241.62 & 61.51 \\
        & & 2 & 300 & - & 183 & 131\,247.19 & 64\,529.69 & 50.83 & 136 & 101\,878.19 & 64\,770.79 & 36.42 \\
        & & 3 & 300 & - & 177 & 127\,266.82 & 64\,495.92 & 49.32 & 132 & 99\,666.82 & 64\,732.82 & 35.05 \\
        & & 4 & 300 & - & 183 & 131\,248.39 & 64\,472.79 & 50.88 & 130 & 98\,710.90 & 64\,720.29 & 34.43 \\
        & & 10 & 300 & - & 176 & 126\,819.47 & 64\,490.87 & 49.15 & 146 & 107\,842.47 & 64\,710.17 & 40.00 \\
        & & 2 & 600 & - & 197 & 139\,868.28 & 64\,810.18 & 53.66 & 140 & 104\,787.28 & 65\,128.28 & 37.85 \\
        & & 3 & 600 & - & 396 & 264\,587.52 & 64\,799.42 & 75.51 & 134 & 100\,769.52 & 65\,080.32 & 35.42 \\
        & & 4 & 600 & - & 392 & 262\,115.66 & 64\,786.06 & 75.28 & 128 & 96\,955.66 & 65\,109.16 & 32.85 \\
        & & 10 & 600 & - & 411 & 273\,771.22 & 64\,610.32 & 76.40 & 128 & 97\,202.92 & 65\,081.62 & 33.05 \\
        \cmidrule{2-13}
        & \multirow{12}{*}{68\,374.20} & 2 & 60 & $\times$ & - & - & - & - & 424 & 281\,642.71 & 64\,107.11 & 77.24 \\
        & & 3 & 60 & $\times$ & - & - & - & - & 407 & 271\,128.64 & 64\,076.44 & 76.37 \\
        & & 4 & 60 & $\times$ & - & - & - & - & 415 & 276\,274.08 & 64\,067.07 & 76.81 \\
        & & 10 & 60 & $\times$ & - & - & - & - & 411 & 273\,573.65 & 64\,046.05 & 76.59 \\
        & & 2 & 300 & $\times$ & - & - & - & - & 136 & 101\,971.64 & 64\,745.54 & 36.51 \\
        & & 3 & 300 & $\times$ & - & - & - & - & 141 & 105\,122.22 & 64\,730.22 & 38.42 \\
        & & 4 & 300 & $\times$ & - & - & - & - & 135 & 101\,668.97 & 64\,700.87 & 36.36 \\
        & & 10 & 300 & $\times$ & - & - & - & - & 146 & 108\,084.75 & 64\,669.95 & 40.17 \\
        & & 2 & 600 & $\times$ & - & - & - & - & 140 & 104\,787.28 & 65\,124.58 & 37.85 \\
        & & 3 & 600 & $\times$ & - & - & - & - & 134 & 100\,769.52 & 65\,078.02 & 35.42 \\
        & & 4 & 600 & $\times$ & - & - & - & - & 128 & 96\,955.66 & 65\,098.56 & 32.86 \\
        & & 10 & 600 & $\times$ & - & - & - & - & 132 & 99\,601.10 & 65\,073.20 & 34.67 \\
        \bottomrule
    \end{tabularx}
    \caption{\Cref{tab:discParam5} continued.}
    \label{tab:discParam6}
\end{table}

\begin{table}[H]
    \scriptsize
    \begin{tabularx}{\textwidth}{c X r r | c r r r r r r r r}
        \toprule
        name & start sol & $\numLinSegments$ & $\timestep$ & fs.? & $\#\course$ & obj. & bound & \makecell[c]{gap\\(\%)} & $\#\course$ & obj. & bound & \makecell[c]{gap\\(\%)} \\
        \midrule
        \multirow{24}{*}{L} & \multirow{12}{*}{-} & 2 & 60 & - & - & - & 80\,040.03 & - & - & - & 80\,044.62 & - \\
        & & 3 & 60 & - & - & - & 80\,041.63 & - & - & - & 80\,045.92 & - \\
        & & 4 & 60 & - & - & - & 80\,042.33 & - & - & - & 80\,045.98 & - \\
        & & 10 & 60 & - & - & - & 80\,009.53 & - & - & - & 80\,046.11 & - \\
        & & 2 & 300 & - & - & - & 80\,043.83 & - & - & - & 80\,044.08 & - \\
        & & 3 & 300 & - & - & - & 80\,043.43 & - & - & - & 80\,043.86 & - \\
        & & 4 & 300 & - & - & - & 80\,045.23 & - & - & - & 80\,045.30 & - \\
        & & 10 & 300 & - & - & - & 80\,044.93 & - & - & - & 80\,045.59 & - \\
        & & 2 & 600 & - & - & - & 80\,043.73 & - & - & - & 80\,043.93 & - \\
        & & 3 & 600 & - & - & - & 80\,044.63 & - & - & - & 80\,044.69 & - \\
        & & 4 & 600 & - & - & - & 80\,045.83 & - & - & - & 80\,045.83 & - \\
        & & 10 & 600 & - & - & - & 80\,044.13 & - & - & - & 80\,044.19 & - \\
        \cmidrule{2-13}
        & \multirow{12}{*}{80\,270.80} & 2 & 60 & $\checkmark$ & 74 & 80\,270.80 & 80\,052.26 & 0.27 & 74 & 80\,096.56 & 80\,058.66 & 0.05 \\
        & & 3 & 60 & $\checkmark$ & 74 & 80\,270.80 & 80\,041.51 & 0.29 & 74 & 80\,270.11 & 80\,047.31 & 0.28 \\
        & & 4 & 60 & $\checkmark$ & 74 & 80\,270.53 & 80\,040.03 & 0.29 & 74 & 80\,250.13 & 80\,045.53 & 0.25 \\
        & & 10 & 60 & $\checkmark$ & 74 & 80\,270.80 & 80\,041.40 & 0.29 & 74 & 80\,264.90 & 80\,046.70 & 0.27 \\
        & & 2 & 300 & $\checkmark$ & 74 & 80\,270.80 & 80\,048.93 & 0.28 & 74 & 80\,253.63 & 80\,049.13 & 0.25 \\
        & & 3 & 300 & $\checkmark$ & 74 & 80\,270.80 & 80\,050.44 & 0.27 & 74 & 80\,126.74 & 80\,050.44 & 0.10 \\
        & & 4 & 300 & $\checkmark$ & 74 & 80\,270.80 & 80\,047.46 & 0.28 & 74 & 80\,106.26 & 80\,047.86 & 0.07 \\
        & & 10 & 300 & $\checkmark$ & 74 & 80\,268.55 & 80\,047.95 & 0.27 & 74 & 80\,142.15 & 80\,050.15 & 0.11 \\
        & & 2 & 600 & $\times$ & 74 & 80\,228.48 & 80\,050.38 & 0.22 & 74 & 80\,173.08 & 80\,050.48 & 0.15 \\
        & & 3 & 600 & $\times$ & 74 & 80\,264.95 & 80\,054.35 & 0.26 & 74 & 80\,113.25 & 80\,054.55 & 0.07 \\
        & & 4 & 600 & $\times$ & 74 & 80\,264.01 & 80\,048.01 & 0.27 & 74 & 80\,132.51 & 80\,048.11 & 0.11 \\
        & & 10 & 600 & $\times$ & 74 & 80\,263.21 & 80\,049.81 & 0.27 & 74 & 80\,110.61 & 80\,050.81 & 0.07 \\
        \midrule
        \multirow{24}{*}{M} & \multirow{12}{*}{-} & 2 & 60 & - & - & - & 61\,399.23 & - & 95 & 76\,809.93 & 61\,426.13 & 20.03 \\
        & & 3 & 60 & - & 160 & 116\,771.76 & 61\,413.66 & 47.41 & 95 & 77\,128.36 & 61\,428.26 & 20.36 \\
        & & 4 & 60 & - & 174 & 125\,444.89 & 61\,400.89 & 51.05 & 108 & 85\,728.29 & 61\,414.59 & 28.36 \\
        & & 10 & 60 & - & 165 & 119\,813.35 & 61\,406.85 & 48.75 & 110 & 85\,862.45 & 61\,421.04 & 28.47 \\
        & & 2 & 300 & - & 132 & 99\,303.73 & 62\,595.23 & 36.97 & 80 & 68\,261.73 & 62\,615.73 & 8.27 \\
        & & 3 & 300 & - & 393 & 262\,970.25 & 62\,596.05 & 76.20 & 84 & 69\,904.15 & 62\,614.55 & 10.43 \\
        & & 4 & 300 & - & 399 & 266\,849.44 & 62\,594.44 & 76.54 & 78 & 66\,909.34 & 62\,625.34 & 6.40 \\
        & & 10 & 300 & - & 371 & 249\,079.95 & 62\,587.05 & 74.87 & 117 & 90\,520.95 & 62\,603.85 & 30.84 \\
        & & 2 & 600 & - & 126 & 96\,304.80 & 63\,221.60 & 34.35 & 78 & 66\,348.90 & 63\,236.60 & 4.69 \\
        & & 3 & 600 & - & 119 & 92\,026.45 & 63\,228.55 & 31.29 & 76 & 65\,426.95 & 63\,237.35 & 3.35 \\
        & & 4 & 600 & - & 392 & 261\,890.58 & 63\,213.88 & 75.86 & 77 & 65\,792.88 & 63\,238.48 & 3.88 \\
        & & 10 & 600 & - & 142 & 106\,274.88 & 63\,222.68 & 40.51 & 86 & 71\,246.68 & 63\,231.08 & 11.25 \\
        \cmidrule{2-13}
        & \multirow{12}{*}{62\,980.50} & 2 & 60 & $\times$ & - & - & - & - & 95 & 76\,809.93 & 61\,426.13 & 20.03 \\
        & & 3 & 60 & $\times$ & - & - & - & - & 96 & 77\,728.35 & 61\,428.75 & 20.97 \\
        & & 4 & 60 & $\times$ & - & - & - & - & 110 & 86\,975.87 & 61\,413.77 & 29.39 \\
        & & 10 & 60 & $\times$ & - & - & - & - & 110 & 85\,862.47 & 61\,420.58 & 28.47 \\
        & & 2 & 300 & $\times$ & 391 & 261\,584.13 & 62\,595.23 & 76.07 & 80 & 68\,261.73 & 62\,615.73 & 8.27 \\
        & & 3 & 300 & $\times$ & - & - & 62\,583.55 & - & 84 & 69\,904.15 & 62\,614.55 & 10.43 \\
        & & 4 & 300 & $\times$ & - & - & 62\,598.60 & - & 78 & 67\,042.40 & 62\,629.40 & 6.58 \\
        & & 10 & 300 & $\times$ & - & - & 62\,564.72 & - & 128 & 97\,318.82 & 62\,594.82 & 35.68 \\
        & & 2 & 600 & $\times$ & 126 & 96\,299.33 & 63\,216.13 & 34.35 & 79 & 67\,075.63 & 63\,231.13 & 5.73 \\
        & & 3 & 600 & $\times$ & 139 & 104\,070.75 & 63\,224.45 & 39.25 & 76 & 65\,426.95 & 63\,237.35 & 3.35 \\
        & & 4 & 600 & $\times$ & - & - & 63\,213.88 & - & 77 & 65\,792.88 & 63\,238.48 & 3.88 \\
        & & 10 & 600 & $\times$ & 142 & 106\,274.88 & 63\,222.38 & 40.51 & 86 & 71\,246.68 & 63\,231.08 & 11.25 \\
        \bottomrule
    \end{tabularx}
    \caption{\Cref{tab:discParam6} continued.}
    \label{tab:discParam7}
\end{table}

\begin{table}[H]
    \scriptsize
    \begin{tabularx}{\textwidth}{c X r r | c r r r r r r r r}
        \toprule
        name & start sol & $\numLinSegments$ & $\timestep$ & fs.? & $\#\course$ & obj. & bound & \makecell[c]{gap\\(\%)} & $\#\course$ & obj. & bound & \makecell[c]{gap\\(\%)} \\
        \midrule
        \multirow{24}{*}{N} & \multirow{12}{*}{-} & 2 & 60 & - & 52 & 37\,123.62 & 32\,388.22 & 12.76 & 43 & 33\,782.22 & 32\,388.22 & 4.13 \\
        & & 3 & 60 & - & 44 & 34\,030.03 & 32\,366.03 & 4.89 & 44 & 33\,897.23 & 32\,366.03 & 4.52 \\
        & & 4 & 60 & - & 42 & 33\,390.75 & 32\,403.15 & 2.96 & 41 & 32\,658.95 & 32\,403.15 & 0.78 \\
        & & 10 & 60 & - & 41 & 32\,787.88 & 32\,369.48 & 1.28 & 41 & 32\,622.18 & 32\,369.48 & 0.77 \\
        & & 2 & 300 & - & 44 & 34\,292.71 & 32\,392.21 & 5.54 & 43 & 33\,684.11 & 32\,392.21 & 3.84 \\
        & & 3 & 300 & - & 43 & 33\,858.70 & 32\,389.00 & 4.34 & 42 & 32\,997.30 & 32\,389.00 & 1.84 \\
        & & 4 & 300 & - & 43 & 33\,718.53 & 32\,399.73 & 3.91 & 43 & 33\,493.63 & 32\,399.73 & 3.27 \\
        & & 10 & 300 & - & 43 & 33\,757.89 & 32\,384.19 & 4.07 & 43 & 33\,418.09 & 32\,384.19 & 3.09 \\
        & & 2 & 600 & - & 44 & 34\,394.40 & 32\,397.00 & 5.81 & 43 & 33\,549.60 & 32\,397.00 & 3.44 \\
        & & 3 & 600 & - & 44 & 34\,570.71 & 32\,388.11 & 6.31 & 44 & 34\,275.11 & 32\,388.11 & 5.51 \\
        & & 4 & 600 & - & 51 & 37\,447.16 & 32\,389.56 & 13.51 & 45 & 35\,291.96 & 32\,389.56 & 8.22 \\
        & & 10 & 600 & - & 44 & 34\,345.01 & 32\,387.41 & 5.70 & 43 & 33\,753.41 & 32\,387.41 & 4.05 \\
        \cmidrule{2-13}
        & \multirow{12}{*}{35\,566.90} & 2 & 60 & $\times$ & 52 & 37\,123.62 & 32\,388.22 & 12.76 & 43 & 33\,782.22 & 32\,388.22 & 4.13 \\
        & & 3 & 60 & $\times$ & 44 & 34\,030.03 & 32\,366.03 & 4.89 & 44 & 33\,897.23 & 32\,366.03 & 4.52 \\
        & & 4 & 60 & $\times$ & 42 & 33\,390.75 & 32\,403.15 & 2.96 & 41 & 32\,658.95 & 32\,403.15 & 0.78 \\
        & & 10 & 60 & $\times$ & 41 & 32\,789.31 & 32\,370.91 & 1.28 & 41 & 32\,634.11 & 32\,370.91 & 0.81 \\
        & & 2 & 300 & $\times$ & 44 & 34\,292.71 & 32\,392.21 & 5.54 & 43 & 33\,684.11 & 32\,392.21 & 3.84 \\
        & & 3 & 300 & $\times$ & 43 & 33\,858.70 & 32\,389.00 & 4.34 & 42 & 32\,997.30 & 32\,389.00 & 1.84 \\
        & & 4 & 300 & $\times$ & 43 & 33\,718.53 & 32\,399.73 & 3.91 & 43 & 33\,493.63 & 32\,399.73 & 3.27 \\
        & & 10 & 300 & $\times$ & 44 & 34\,194.09 & 32\,384.19 & 5.29 & 43 & 33\,418.09 & 32\,384.19 & 3.09 \\
        & & 2 & 600 & $\times$ & 44 & 34\,394.40 & 32\,397.00 & 5.81 & 43 & 33\,549.60 & 32\,397.00 & 3.44 \\
        & & 3 & 600 & $\times$ & 44 & 34\,570.71 & 32\,388.11 & 6.31 & 44 & 34\,275.11 & 32\,388.11 & 5.51 \\
        & & 4 & 600 & $\times$ & 51 & 37\,447.16 & 32\,389.56 & 13.51 & 45 & 35\,291.96 & 32\,389.56 & 8.22 \\
        & & 10 & 600 & $\times$ & 44 & 34\,478.71 & 32\,387.41 & 6.07 & 43 & 33\,753.41 & 32\,387.41 & 4.05 \\
        \midrule
        \multirow{24}{*}{O} & \multirow{12}{*}{-} & 2 & 60 & - & - & - & 36\,944.20 & - & 34 & 37\,010.40 & 36\,955.50 & 0.15 \\
        & & 3 & 60 & - & - & - & 36\,949.92 & - & 34 & 37\,016.62 & 36\,958.52 & 0.16 \\
        & & 4 & 60 & - & - & - & 36\,945.37 & - & 34 & 37\,012.77 & 36\,954.37 & 0.16 \\
        & & 10 & 60 & - & - & - & 36\,961.37 & - & 34 & 37\,014.57 & 36\,969.17 & 0.12 \\
        & & 2 & 300 & - & - & - & 36\,945.79 & - & 34 & 37\,017.59 & 36\,952.69 & 0.18 \\
        & & 3 & 300 & - & - & - & 36\,937.12 & - & 34 & 37\,017.52 & 36\,944.62 & 0.20 \\
        & & 4 & 300 & - & 34 & 37\,047.81 & 36\,956.51 & 0.25 & 34 & 37\,019.41 & 36\,968.71 & 0.14 \\
        & & 10 & 300 & - & - & - & 36\,945.13 & - & 34 & 37\,011.63 & 36\,953.63 & 0.16 \\
        & & 2 & 600 & - & 34 & 37\,034.01 & 36\,938.51 & 0.26 & 34 & 37\,015.81 & 36\,941.71 & 0.20 \\
        & & 3 & 600 & - & - & - & 36\,943.16 & - & 34 & 37\,018.16 & 36\,949.56 & 0.19 \\
        & & 4 & 600 & - & - & - & 36\,947.55 & - & 34 & 37\,018.15 & 36\,956.45 & 0.17 \\
        & & 10 & 600 & - & 34 & 37\,051.26 & 36\,930.06 & 0.33 & 34 & 37\,018.76 & 36\,937.36 & 0.22 \\
        \cmidrule{2-13}
        & \multirow{12}{*}{37\,014.00} & 2 & 60 & $\checkmark$ & 34 & 37\,014.00 & 36\,968.76 & 0.12 & 34 & 37\,013.96 & 36\,970.56 & 0.12 \\
        & & 3 & 60 & $\checkmark$ & 34 & 37\,014.00 & 36\,975.43 & 0.10 & 34 & 37\,013.43 & 36\,976.83 & 0.10 \\
        & & 4 & 60 & $\checkmark$ & 34 & 37\,013.95 & 36\,976.25 & 0.10 & 34 & 37\,013.95 & 36\,978.05 & 0.10 \\
        & & 10 & 60 & $\checkmark$ & 34 & 37\,014.00 & 36\,978.72 & 0.10 & 34 & 37\,014.00 & 36\,979.82 & 0.09 \\
        & & 2 & 300 & $\checkmark$ & 34 & 37\,014.00 & 36\,968.48 & 0.12 & 34 & 37\,013.98 & 36\,968.48 & 0.12 \\
        & & 3 & 300 & $\checkmark$ & 34 & 37\,014.00 & 36\,976.61 & 0.10 & 34 & 37\,014.00 & 36\,978.21 & 0.10 \\
        & & 4 & 300 & $\checkmark$ & 34 & 37\,014.00 & 36\,976.41 & 0.10 & 34 & 37\,014.00 & 36\,977.21 & 0.10 \\
        & & 10 & 300 & $\checkmark$ & 34 & 37\,013.63 & 36\,977.03 & 0.10 & 34 & 37\,013.43 & 36\,977.73 & 0.10 \\
        & & 2 & 600 & $\times$ & 34 & 37\,051.21 & 36\,938.31 & 0.30 & 34 & 37\,015.81 & 36\,941.71 & 0.20 \\
        & & 3 & 600 & $\checkmark$ & 34 & 37\,014.00 & 36\,969.21 & 0.12 & 34 & 37\,014.00 & 36\,969.21 & 0.12 \\
        & & 4 & 600 & $\checkmark$ & 34 & 37\,014.00 & 36\,970.58 & 0.12 & 34 & 37\,014.00 & 36\,970.58 & 0.12 \\
        & & 10 & 600 & $\checkmark$ & 34 & 37\,014.00 & 36\,976.92 & 0.10 & 34 & 37\,014.00 & 36\,976.92 & 0.10 \\
        \bottomrule
    \end{tabularx}
    \caption{\Cref{tab:discParam7} continued.}
    \label{tab:discParam8}
\end{table}

\begin{table}[H]
    \scriptsize
    \begin{tabularx}{\textwidth}{c X r r | c r r r r r r r r}
        \toprule
        name & start sol & $\numLinSegments$ & $\timestep$ & fs.? & $\#\course$ & obj. & bound & \makecell[c]{gap\\(\%)} & $\#\course$ & obj. & bound & \makecell[c]{gap\\(\%)} \\
        \midrule
        \multirow{24}{*}{P} & \multirow{12}{*}{-} & 2 & 60 & - & - & - & 283\,888.81 & - & - & - & 284\,271.76 & - \\
        & & 3 & 60 & - & - & - & 283\,888.81 & - & - & - & 284\,314.25 & - \\
        & & 4 & 60 & - & - & - & 283\,888.81 & - & - & - & 284\,299.49 & - \\
        & & 10 & 60 & - & - & - & 284\,130.81 & - & - & - & 284\,212.82 & - \\
        & & 2 & 300 & - & - & - & 284\,154.81 & - & - & - & 284\,251.48 & - \\
        & & 3 & 300 & - & - & - & 284\,155.81 & - & - & - & 284\,285.42 & - \\
        & & 4 & 300 & - & - & - & 284\,150.81 & - & - & - & 284\,269.99 & - \\
        & & 10 & 300 & - & - & - & 284\,132.81 & - & - & - & 284\,235.96 & - \\
        & & 2 & 600 & - & - & - & 284\,217.81 & - & - & - & 284\,276.61 & - \\
        & & 3 & 600 & - & - & - & 284\,163.81 & - & - & - & 284\,287.99 & - \\
        & & 4 & 600 & - & - & - & 284\,156.81 & - & - & - & 284\,277.92 & - \\
        & & 10 & 600 & - & - & - & 284\,084.81 & - & - & - & 284\,259.57 & - \\
        \cmidrule{2-13}
        & \multirow{12}{*}{286\,061.00} & 2 & 60 & $\times$ & 83 & 286\,067.30 & 283\,893.30 & 0.76 & 83 & 286\,066.30 & 284\,277.30 & 0.63 \\
        & & 3 & 60 & $\checkmark$ & 83 & 286\,061.00 & 283\,892.44 & 0.76 & 83 & 286\,060.44 & 284\,284.44 & 0.62 \\
        & & 4 & 60 & $\checkmark$ & 83 & 286\,061.00 & 283\,897.80 & 0.76 & 83 & 286\,043.80 & 284\,263.80 & 0.62 \\
        & & 10 & 60 & $\checkmark$ & 83 & 286\,060.98 & 284\,133.98 & 0.67 & 83 & 286\,060.98 & 284\,246.98 & 0.63 \\
        & & 2 & 300 & $\times$ & 83 & 286\,055.06 & 284\,152.06 & 0.67 & 83 & 286\,049.06 & 284\,253.06 & 0.63 \\
        & & 3 & 300 & $\checkmark$ & 83 & 286\,028.83 & 284\,058.83 & 0.69 & 83 & 286\,017.83 & 284\,291.83 & 0.60 \\
        & & 4 & 300 & $\checkmark$ & 83 & 286\,061.00 & 284\,055.89 & 0.70 & 82 & 285\,673.89 & 284\,264.89 & 0.49 \\
        & & 10 & 300 & $\checkmark$ & 83 & 286\,061.00 & 284\,134.42 & 0.67 & 83 & 286\,061.00 & 284\,198.42 & 0.65 \\
        & & 2 & 600 & $\times$ & 83 & 286\,036.71 & 284\,226.71 & 0.63 & 82 & 285\,715.71 & 284\,294.71 & 0.50 \\
        & & 3 & 600 & $\times$ & 82 & 285\,693.76 & 284\,162.76 & 0.54 & 82 & 285\,668.76 & 284\,287.76 & 0.48 \\
        & & 4 & 600 & $\times$ & 82 & 285\,643.21 & 284\,089.21 & 0.54 & 82 & 285\,636.21 & 284\,289.21 & 0.47 \\
        & & 10 & 600 & $\times$ & 83 & 286\,028.81 & 284\,059.81 & 0.69 & 82 & 285\,656.81 & 284\,263.81 & 0.49 \\
        \midrule
        \bottomrule
    \end{tabularx}
    \caption{\Cref{tab:discParam8} continued.}
    \label{tab:discParam9}
\end{table}

\end{document}